\documentclass[reqno]{amsart}
\usepackage{amsmath, amsthm, amscd, amssymb, amsfonts, amsbsy}
\usepackage{mathrsfs}
\usepackage{pxfonts}
\usepackage{bbm}
\usepackage{todonotes}
\usepackage{marginnote}

\theoremstyle{plain}
\newtheorem{theorem}[equation]{Theorem}
\newtheorem{lemma}[equation]{Lemma}
\newtheorem{corollary}[equation]{Corollary}
\newtheorem{proposition}[equation]{Proposition}

\theoremstyle{definition}
\newtheorem{definition}[equation]{Definition}
\newtheorem{condition}[equation]{Condition}

\theoremstyle{remark}

\newcommand{\capacity}{\operatorname{cap}}

\newcommand{\supp}{\operatorname{supp}}

\numberwithin{equation}{section}

\newcommand\cD{\mathcal{D}}
\newcommand\cE{\mathcal{E}}

\providecommand{\abs}[1]{\lvert#1\rvert}

\providecommand{\norm}[1]{\lVert#1\rVert}

\newcommand{\red}{\mathscr R}
\newcommand{\bal}{\hat{\mathscr R}}

\begin{document}

\title[Wiener-type criteria for parabolic equations]{Wiener-type criteria for parabolic equations with two-sided Gaussian bounds}

\author[P. Jung]{Pilgyu Jung}
\address[P. Jung]{Department of Mathematics, Yonsei University, 50 Yonsei-Ro, Seodaemun-gu, Seoul 03722, Republic of Korea}
\email{pilgyujung@yonsei.ac.kr}

\author[S. Kim]{Seick Kim}
\address[S. Kim]{Department of Mathematics, Yonsei University, 50 Yonsei-Ro, Seodaemun-gu, Seoul 03722, Republic of Korea}
\email{kimseick@yonsei.ac.kr}
\thanks{S. Kim is supported by the National Research Foundation of Korea (NRF) under agreement NRF-2022R1A2C1003322.
S. Kim is also affiliated with the Intelligent Computational Science Institute (IN2CSI) at Yonsei University.}

\author[S. Lee]{Sungjin Lee}
\address[S. Lee]{Department of Mathematics, Sogang University, 35 Baekbeom-ro, Mapo-gu, Seoul 04107, Republic of Korea}
\email{sungjinlee@sogang.ac.kr}

\date{}


\begin{abstract}
We establish sufficient and necessary Wiener-type criteria for boundary regularity of uniformly parabolic equations whose fundamental solutions satisfy two-sided Gaussian bounds.
In particular, our sufficient criterion holds at the critical exponent in the Gaussian lower bound.
For the heat equation, the resulting double-series criterion is equivalent to the classical Wiener criterion of Evans and Gariepy.
Our results apply to divergence form equations with bounded measurable, possibly nonsymmetric coefficients, and to non-divergence form equations with coefficients having Dini mean oscillation in the spatial variables and merely measurable in time.
As an application, we characterize boundary regularity in a class of time-varying domains in terms of elliptic boundary regularity of their spatial cross-sections.
We also establish a quantitative boundary estimate in terms of finite partial sums of the Wiener-type series at the critical exponent.
\end{abstract}

\maketitle

\section{Introduction and main results}
We consider a parabolic operator $P$ on $\mathbb R^{d+1}$ of the form
\[
Pu=\partial_t u-\sum_{i,j=1}^d D_i(a^{ij}(t,x)D_ju)
\quad\text{ or }\quad
Pu=\partial_t u-\sum_{i,j=1}^d a^{ij}(t,x)D_{ij}u,
\]
in divergence or non-divergence form, respectively.
In either case, the coefficient matrix $\mathbf A=(a^{ij})$ is real-valued and measurable, and satisfies the uniform ellipticity and boundedness conditions
\begin{equation}	\label{eq0208sat}
\mathbf A(t,x)\xi\cdot\xi\ge\lambda\abs{\xi}^2,  \qquad \abs{\mathbf A(t,x)\xi} \le \Lambda \abs{\xi},  \qquad 0<\lambda \le \Lambda<\infty,
\end{equation}
for all $\xi \in \mathbb{R}^d$ and $(t,x) \in \mathbb{R}^{d+1}$.
In the non-divergence case, we assume without loss of generality that  $\mathbf A$ is symmetric.

A fundamental solution $\Gamma(t,x,\tau,\xi)$ for $P$ is a function that vanishes for $t\le\tau$ and formally satisfies
\[
P\Gamma(\cdot,\cdot,\tau,\xi)=\delta_{(\tau,\xi)} \quad\text{in }\mathbb R^{d+1}.
\]
Equivalently, it satisfies
\[
P\Gamma(\cdot,\cdot,\tau,\xi)=0 \quad\text{in }(\tau,\infty)\times\mathbb R^d,
\]
and
\[
\lim_{t\downarrow\tau}\Gamma(t,\cdot,\tau,\xi)=\delta_\xi \quad\text{on }\mathbb R^d.
\]
We further assume that $\Gamma$ satisfies the \emph{mass conservation identity}
\[
\int_{\mathbb R^d}\Gamma(t,x,\tau,\xi)\,d\xi=1,\qquad t>\tau,
\]
and the \emph{reproduction identity}
\[
\Gamma(t,x,\tau,\xi)=\int_{\mathbb R^d} \Gamma(t,x,s,y)\Gamma(s,y, \tau,\xi)\,dy,\qquad \tau<s<t.
\]

For $a>0$, define the Gaussian kernel
\[
G_a(t,x)=
\begin{cases}
t^{-d/2}\exp(-a\abs{x}^2/t),&t>0,\\ 0,&t \le 0.
\end{cases}
\]
For $P$ in divergence form, Aronson's estimates \cite{Aronson} yield a fundamental solution satisfying the two-sided Gaussian bounds
\begin{equation}	\label{eq0157sat}
c_0 G_{b_0}(t-\tau, x-\xi) \le \Gamma(t,x,\tau,\xi) \le C_0 G_{a_0}(t-\tau, x-\xi),
\end{equation}
where $0<a_0\le b_0$ and $c_0,C_0>0$ depend only on $d$, $\lambda$, and $\Lambda$.
These bounds hold for bounded measurable coefficients, without any symmetry assumption on $\mathbf A$.

In contrast, for $P$ in non-divergence form, two-sided Gaussian bounds cannot in general be expected for merely measurable coefficients.
Indeed, such bounds may fail even when the coefficients are continuous in both $t$ and $x$, as demonstrated by the counterexample in \cite{Ilin62}; see also \cite{FK81, Sa81}.

For $P$ in non-divergence form with $\mathbf A\in \mathrm{DMO}_x$ (see Definition~\ref{dmo_x}), two-sided Gaussian bounds for the fundamental solution are given in Remark~5.12 of \cite{DEK21}.
These follow by combining the adjoint estimates of \cite{DEK21} with the Gaussian estimates of Escauriaza~\cite{Esc2000}; see also \cite{KLS}.
Such bounds were previously established under H\"older or Dini continuity assumptions on the coefficients; see, for instance, \cite{Friedman,PE84}.

The Wiener criterion for parabolic equations has a long history, with early contributions by Pini \cite{Pini}.
Landis \cite{Landis} obtained a necessary and sufficient criterion for the heat equation in terms of caloric equilibrium potentials.
Lanconelli \cite{Lanconelli1973} established a necessary condition involving thermal capacity, and Evans and Gariepy \cite{EG1982} subsequently proved the corresponding necessary and sufficient criterion.
Extensions to parabolic operators with smooth variable coefficients were obtained by Garofalo and Lanconelli \cite{GL1988}, and to divergence form operators with $C^1$-Dini continuous coefficients by Fabes, Garofalo, and Lanconelli \cite{FGL1989}.

For divergence form parabolic equations with measurable coefficients, Lanconelli \cite{Lanconelli1975a} established a Gaussian capacity double series criterion characterizing boundary points regular for every operator in the class.
In a companion paper \cite{Lanconelli1975b}, he also obtained quantitative estimates for the boundary modulus of continuity of solutions to the Dirichlet problem.
Lanconelli \cite{Lanconelli1977} subsequently obtained anisotropic Wiener-type criteria for non-divergence form operators under strict inequalities between the associated quadratic forms.
These works laid the groundwork for subsequent Gaussian capacity approaches to parabolic Wiener tests.

Lanconelli, Tralli, and Uguzzoni \cite[Theorem~1.1]{LTU} extended this theory to a class of hypoelliptic diffusion operators satisfying two-sided Gaussian bounds.
In the Euclidean setting, Gaussian capacities associated with different positive exponents are equivalent.
Their sufficient criterion requires $b>b_0$, whereas the necessary criterion holds for $a\le a_0$.
In the present paper, we establish sufficiency at the critical exponent $b=b_0$, thereby obtaining a necessary and sufficient criterion when $a_0=b_0$.
In particular, for the heat equation, the resulting Lanconelli-type double series criterion is equivalent to the classical Evans--Gariepy Wiener criterion.
We also establish a quantitative boundary estimate in terms of finite partial sums of the Wiener-type series at the critical exponent $b_0$, extending the corresponding estimate in \cite[Theorem~1.2]{LTU} to the endpoint in the Euclidean setting.

\subsection{Two-sided Gaussian bounds}
We impose the following condition.

\begin{condition}	\label{gaussian}
The operator $P$ admits a fundamental solution $\Gamma(t,x,\tau,\xi)$ on $\mathbb R^{d+1}\times\mathbb R^{d+1}$ such that, for some $T_0\in(0,\infty]$, $0<a_0\le b_0$, and $c_0,C_0>0$, the two-sided Gaussian bounds \eqref{eq0157sat} hold whenever $0<t-\tau<T_0$.
\end{condition}

By the reproduction identity, the Gaussian bounds hold for
$0<t-\tau<T$ for every finite $T>0$, with the same exponents $a_0 ,b_0$ but possibly different constants $c_0, C_0$, depending only on $d$, $T$, and the constants in Condition~\ref{gaussian}.
Throughout the paper, we use this fact without further comment.

For bounded continuous $\psi$ and bounded compactly supported $f$, define
\[
\mathscr S_\tau\psi(t,x)=\int_{\mathbb R^d}\Gamma(t,x,\tau,\xi)\psi(\xi)\,d\xi, \qquad t>\tau,
\]
and
\[
\mathscr Gf(t,x)=\int_{-\infty}^t\int_{\mathbb R^d}\Gamma(t,x,\tau,\xi)f(\tau,\xi)\,d\xi\,d\tau.
\]
By standard parabolic theory, $\mathscr S_\tau\psi$ solves $Pu=0$ with initial datum $\psi$ at $t=\tau$, while $\mathscr Gf$ is a bounded solution of $Pu=f$.

\subsection{Gaussian capacities and LTU sets}
For $a>0$ and a nonnegative Borel measure $\mu$, define
\[
G_a\mu(z)=\int_{\mathbb R^{d+1}}
G_a(z-\zeta)\,d\mu(\zeta).
\]
For a compact set $E\subset\mathbb R^{d+1}$, define
\[
\mathcal C_a(E)=\sup\left\{\mu(E):\
\mu\ge0,\; \supp\mu\subset E,\;
G_a\mu\le1\text{ on }\mathbb R^{d+1}\right\}.
\]
The classical thermal capacity is given by
\[
\capacity(E):=(4\pi)^{d/2}\mathcal C_{1/4}(E),
\]
corresponding to the standard heat kernel.
As shown in \cite{Lanconelli1975a},
\[
\mathcal C_{a}(E)\asymp\capacity(E), \qquad a>0,
\]
where the comparison constants depend only on $a$ and $d$.
Throughout, $A\asymp B$ means that
\[
cA\le B\le CA
\]
for some constants $c,C>0$ independent of the quantities being compared.

For $z=(t,x)\in\mathbb R^{d+1}$, set
\[
\abs{z}=\max\{\abs{x},\sqrt{\abs{t}}\},
\]
so that $\abs{z-\zeta}$ denotes the parabolic distance between $z$ and $\zeta$.

For $z_0=(t_0,x_0)\in\mathbb R^{d+1}$ and $r>0$, define
\[
\begin{aligned}
Q_r(z_0)&=\{z:\abs{z-z_0}<r\}=(t_0-r^2,t_0+r^2)\times B_r(x_0),\\
Q_r^-(z_0)&=(t_0-r^2,t_0)\times B_r(x_0).
\end{aligned}
\]

Let $\Omega\subset\mathbb R^{d+1}$ be an open set, and denote its complement by $\Omega^\complement$.
For $z_0=(t_0,x_0)\in\partial\Omega$, $0<\rho<1$, and $k,\ell\in\mathbb N$, define the following variants of the LTU sets:
\begin{multline}	\label{eq0237sat}
\Omega_k^\ell(z_0,\rho)=\biggl\{(t,x)\in\Omega^\complement\cap\overline Q_{\sqrt\rho}(z_0): \rho^{k+1}\le t_0-t\le\rho^k,\\
(\ell-1)\log(1/\rho)\le \frac{\abs{x-x_0}^2}{t_0-t} \le\ell\log(1/\rho)\biggr\}.
\end{multline}
Here $k$ indexes time slabs at geometric scales, while $\ell$ indexes the regions between the backward paraboloids
\[
\abs{x-x_0}^2=(\ell-1)(t_0-t)\log(1/\rho) \quad\text{ and }\quad \abs{x-x_0}^2=\ell(t_0-t)\log(1/\rho).
\]
In particular, $\ell$ indexes the exponential spatial factor through
\[
\rho^\ell\le\exp\left(-\frac{\abs{x-x_0}^2}{t_0-t}\right) \le\rho^{\ell-1}.
\]
Consequently, for any fixed $a>0$,
\[
G_a(t_0-t,x_0-x)\asymp\rho^{a\ell-kd/2} \qquad\text{on }\Omega_k^\ell(z_0,\rho),
\]
with comparison constants depending only on $d$, $a$, and $\rho$.

The double-index decomposition into geometric time slabs and Gaussian spatial levels goes back to Lanconelli \cite{Lanconelli1975a}.
The LTU sets \cite[(1.6)]{LTU} incorporate an additional parabolic distance cutoff.
Our sets differ from theirs only in the choice of parabolic distance, which does not affect the divergence of the corresponding double series.
\subsection{Main results and applications}
\begin{definition}[$\mathrm{DMO}_x$ condition]	\label{dmo_x}
We say that $\mathbf A$ has Dini mean oscillation in $x$, written $\mathbf A\in\mathrm{DMO}_x$, if
\[
\int_0^1\frac{\omega_{\mathbf A}^x(r)}r\,d r<\infty,
\]
where
\[
\omega_{\mathbf A}^x(r)=\sup_{z=(t,x)\in\mathbb R^{d+1}} \fint_{t-r^2}^{t} \fint_{B_r(x)} \abs{\mathbf A(\tau,\xi)-\bar{\mathbf A}_{x,r}(\tau)} \,d\xi\,d\tau,\qquad
\bar{\mathbf A}_{x,r}(\tau)=\fint_{B_r(x)} \mathbf A(\tau,\xi)\,d\xi.
\]
\end{definition}

Uniform Dini continuity of $\mathbf A$ in $x$ implies $\mathbf A\in\mathrm{DMO}_x$.
While no regularity in $t$ can be expected from the $\mathrm{DMO}_x$ condition, it is worth noting that the condition does not guarantee continuity in $x$ at each fixed time either.

In the following theorem, boundary regularity is understood in the Perron sense; see Section~\ref{sec2} for its precise definition.

\begin{theorem} 	\label{thm01}
Let $P$ be an operator in divergence form with measurable, possibly nonsymmetric $\mathbf A$ satisfying \eqref{eq0208sat}, or in non-divergence form with symmetric $\mathbf A\in\mathrm{DMO}_x$ satisfying \eqref{eq0208sat}.
Then Condition~\ref{gaussian} holds.

Moreover, for any bounded open set $\Omega\subset\mathbb R^{d+1}$, $z_0\in\partial\Omega$, and $0<\rho<1$, the following statements hold.
\begin{enumerate}
\item[(i)]
If
\begin{equation}	\label{eq:sufficient}
\sum_{k=1}^\infty\sum_{\ell=1}^\infty \capacity(\Omega_k^\ell(z_0,\rho)) \rho^{b_0\ell - kd/2}=\infty,
\end{equation}
then $z_0$ is regular.

\item[(ii)]
If $z_0$ is regular, then
\begin{equation}	\label{eq:necessary}
\sum_{k=1}^\infty\sum_{\ell=1}^\infty \capacity(\Omega_k^\ell(z_0,\rho)) \rho^{a_0\ell - kd/2}=\infty.
\end{equation}
\end{enumerate}
\end{theorem}

Our main theorem has elliptic counterparts for divergence form operators with bounded measurable coefficients \cite{LSW63, GW82} and non-divergence form operators with Dini mean oscillation coefficients \cite{DKK25a,DKK25b}; see also \cite{Bauman85}.
An important distinction is that, unlike in the elliptic setting, parabolic boundary regularity may depend on the diffusion coefficient even for constant coefficient heat operators; see \cite{Petrowsky, EK1970, Lanconelli1977}.
Thus, the Gaussian exponents play an essential role, and the sufficient and necessary conditions in Theorem~\ref{thm01} involve the possibly different exponents $b_0$ and $a_0$, respectively.
When $a_0=b_0$, these conditions give a necessary and sufficient criterion.
In particular, we have the following corollary for the heat operator.

\begin{corollary}	\label{cor01}
Let $\Omega\subset\mathbb R^{d+1}$ be bounded and open, and let $z_0\in\partial\Omega$.
For the heat operator $P=\partial_t-\Delta$, the point $z_0$ is regular if and only if
\[
\sum_{k=1}^\infty\sum_{\ell=1}^\infty \capacity(\Omega_k^\ell(z_0,\rho)) \rho^{\ell/4-kd/2}=\infty
\]
for some (equivalently, every) $0<\rho<1$.
\end{corollary}

For comparison, let $z_0=(t_0,x_0)$ and define
\[
A_k(z_0)=\{(t,x):t<t_0,\; 2^{kd/2}\le G_{1/4}(t_0-t,x_0-x) \le 2^{(k+1)d/2}\}.
\]
The Evans--Gariepy criterion \cite{EG1982} states that $z_0$ is regular if and only if
\[
\sum_{k=1}^\infty 2^{kd/2} \capacity(\Omega^\complement\cap \overline A_k(z_0))=\infty.
\]
Thus, the double series criterion in Corollary~\ref{cor01} is equivalent to the Evans--Gariepy single series criterion.
Since the LTU sets are independent of the Gaussian exponents, their formulation is particularly convenient for operators satisfying two-sided Gaussian bounds.

As another application of Theorem~\ref{thm01}, we characterize boundary regularity in certain time-varying domains in terms of elliptic boundary regularity of their spatial cross-sections.
The domains under consideration are locally obtained from cylindrical domains by time-dependent changes of spatial variables.
Let $I\subset\mathbb R$ be an open interval and $B\subset\mathbb R^d$ be a ball.
We consider maps $\Phi:I\times B\to\mathbb R^d$ satisfying
\begin{equation}	\label{eq1254mon}
\begin{aligned}
N^{-1}\abs{x-y}\le \abs{\Phi(t,x)-\Phi(t,y)} &\le N\abs{x-y},\\
\abs{\Phi(t,x)-\Phi(s,x)}&\le N\abs{t-s}^{1/2}
\end{aligned}
\end{equation}
for some $N\ge 1$ and all $t,s\in I$ and $x,y\in B$.

\begin{corollary}	\label{cor02}
Let $P$ satisfy the assumptions of Theorem~\ref{thm01}, and let $\Omega\subset\mathbb R^{d+1}$ be open.
Suppose $z_0=(t_0,x_0)\in\partial\Omega$ and $y_0\in\partial D$, where $D\subset\mathbb R^d$ is open.
Assume that for some $r>0$, a map $\Phi$ defined on
\[
U=Q_r(t_0,y_0)
\]
satisfies \eqref{eq1254mon} with a constant $N \ge 1$.
Set $\mathsf T(t,y)=(t,\Phi(t,y))$ and suppose that $\mathsf T(t_0,y_0)=z_0$ and, for some $\delta>0$ with $Q_\delta^-(z_0)\subset \mathsf T(U)$,
\[
\Omega\cap Q_\delta^-(z_0) =\mathsf T(U\cap(\mathbb R\times D)) \cap Q_\delta^-(z_0).
\]
Then $z_0$ is regular for $P$ in $\Omega$ if and only if $y_0$ is regular for the Laplacian in $D$.
In particular, $z_0$ is regular for $P$ if and only if it is regular for the heat operator $\partial_t-\Delta$.
\end{corollary}

The preliminary results are developed in Sections~\ref{sec2} and~\ref{sec3}.
Theorem~\ref{thm01} and Corollary~\ref{cor02} are proved in Section~\ref{sec4}.
In Section~\ref{sec05}, we establish the quantitative boundary estimate by adapting the technique developed in \cite{Kim2026} for the heat equation.

\section{Parabolic potential theory}	\label{sec2}
We consider weak solutions for operators in divergence form and strong solutions in $W^{1,2}_{d+1,\mathrm{loc}}$ for operators in non-divergence form.
By the De Giorgi--Nash--Moser and Krylov--Safonov theories, respectively, these solutions admit locally H\"older continuous representatives, which we use throughout.
Nonnegative solutions satisfy interior Harnack inequalities.
We also use the corresponding strong minimum principles for nonnegative supersolutions.

For completeness, proofs of several standard lemmas are given in Appendix~\ref{sec_a}.

\subsection{Superparabolic functions and Perron solutions}

A solution of $Pu=0$ is called $P$-caloric.
For a cylinder $Q=(s,t)\times B$, where $B$ is a ball, define its parabolic boundary by
\[
\partial_p Q=\bigl(\{s\}\times\overline B\bigr)\cup\bigl([s,t] \times\partial B\bigr).
\]

For every cylinder $Q$, the Dirichlet problem for $P$ with bounded continuous data on $\partial_pQ$ has a unique bounded $P$-caloric solution.
There exists, for each $z\in Q$, a probability measure $\omega_Q^z$ on $\partial_pQ$, called the $P$-caloric measure, such that the solution with boundary data $\varphi$ is given by
\[
u(z)=\int_{\partial_pQ}\varphi(\zeta)\, d\omega_Q^z(\zeta).
\]

Let $D\subset\mathbb R^{d+1}$ be open.
A function $v:D\to(-\infty,\infty]$ is called $P$-superparabolic if it is lower semicontinuous, finite on a dense subset of $D$, and, for every cylinder $Q$ compactly contained in $D$ and every $P$-caloric function $u\in C(\overline Q)$,
\[
u\le v\quad\text{on }\;\partial_p Q \implies u\le v\quad\text{in }\,Q.
\]

For a finite-valued $P$-superparabolic function $v$ in an open set $D$ and a cylinder $Q\Subset D$, set
\[
\mathscr H_Q v(z)=\int_{\partial_p Q}v(\zeta)\,d\omega_Q^z(\zeta),\qquad z\in Q.
\]
Approximating $v$ on $\partial_pQ$ from below by an increasing sequence of continuous functions, comparison and monotone convergence give $\mathscr H_Qv \le v$ in $Q$.
In particular, the integral is finite.
The Poisson modification of $v$ in $Q$ is defined by
\[
v^Q(z)=
\begin{cases}
\mathscr H_Q v(z),&z\in(s,t]\times B,\\[4pt]
v(z),&\text{otherwise},
\end{cases}
\]
where $\mathscr H_Qv$ on $\{t\}\times B$ is
understood by continuous extension from $Q$.

\begin{lemma}	\label{lem03}
Let $D\subset\mathbb R^{d+1}$ be open.
\begin{enumerate}
\item[(a)]
(Comparison principle)
Suppose $D$ is bounded and $v$ is bounded below and $P$-superparabolic in $D$.
If
\[
\liminf_{D\ni z\to\zeta}v(z)\ge 0 \qquad(\zeta\in\partial D),
\]
then $v \ge 0$ in $D$.
Moreover, if $D\subset\{t<T\}$, it suffices to impose the boundary condition only on $\partial D\cap\{t<T\}$.

\item[(b)]
If $u,v$ are finite-valued $P$-superparabolic functions in $D$ and $c>0$, then $cu$, $u+v$, and $\min\{u,v\}$ are $P$-superparabolic in $D$.

\item[(c)]
If $(v_n)$ is an increasing sequence of $P$-superparabolic functions in $D$ converging pointwise to a finite-valued function $v$, then $v$ is $P$-superparabolic in $D$.

\item[(d)]
(Poisson modification)
If $v$ is finite-valued and $P$-superparabolic in $D$, and $Q\Subset D$ is a cylinder, then $v^Q$ is $P$-superparabolic in $D$, $P$-caloric in $Q$, and satisfies $v^Q\le v$ in $D$.
\end{enumerate}
\end{lemma}

In applications of the comparison principle in a cylinder $Q=(s,t)\times B$, the boundary condition on the terminal edge $\{t\}\times\partial B$ may be omitted by part (a).

\begin{lemma}	\label{lem04}
Let $D\subset\mathbb R^{d+1}$ be open, and let $\mathcal F$ be a nonempty family of finite-valued $P$-superparabolic functions in $D$, uniformly bounded below and closed under finite minima and Poisson modifications in cylinders $Q\Subset D$.
Then $\inf_{v\in\mathcal F}v$ is $P$-caloric in $D$.
\end{lemma}

Let $\Omega\subset\mathbb R^{d+1}$ be bounded and open, and let $f\in C(\partial\Omega)$.
The \emph{upper Perron class} $\mathcal U_f$ consists of all $P$-superparabolic functions $v$ in $\Omega$ that are bounded below and satisfy
\[
\liminf_{\Omega\ni z\to\zeta}v(z)\ge f(\zeta) \qquad(\zeta\in\partial\Omega).
\]
The upper Perron solution is defined by
\[
\overline H_f=\inf_{v\in\mathcal U_f} v.
\]
Here, the infimum may be taken over bounded functions in $\mathcal U_f$.
The lower Perron solution is defined by
\[
\underline H_f=-\overline H_{-f}.
\]
By the comparison principle and Lemma~\ref{lem04}, both $\overline H_f$ and $\underline H_f$ are $P$-caloric in $\Omega$ and satisfy
\[
\min_{\partial\Omega}f \le\underline H_f\le\overline H_f \le\max_{\partial\Omega}f.
\]
A point $z_0\in\partial\Omega$ is called \emph{regular} if
\[
\lim_{\Omega\ni z\to z_0}\overline H_f(z)=f(z_0)
\]
for every $f\in C(\partial\Omega)$.

\subsection{Capacities and reduced functions}
For a compact set $K\subset\mathbb R^{d+1}$, define its \emph{reduced function} by
\[
\red_K(z) =\inf\{v(z):v \ge 0\text{ is $P$-superparabolic in }\mathbb R^{d+1},\ v\ge 1\text{ on }K\}.
\]
For an arbitrary set $E\subset\mathbb R^{d+1}$, define
\[
\red_E(z) =\sup\{\red_K(z):K\subset E,\ K\text{ compact}\}.
\]

\begin{lemma}	\label{lem01}
For any $a>0$ and compact $E\subset\mathbb R^{d+1}$,
\[
\mathcal C_a(E)\asymp\capacity(E)\asymp\mathcal C_a^\dagger(E),
\]
where $\mathcal C_a^\dagger$ denotes the capacity associated
with the adjoint Gaussian kernel
\begin{equation}\label{eq0205sat}
G_a^\dagger(t,x)=G_a(-t,-x)=G_a(-t,x).
\end{equation}
The comparison constants depend only on $d$ and $a$.

Thermal capacity is monotone, countably subadditive under compact coverings, and continuous from above on decreasing sequences of compact sets.
For $r>0$, write
\begin{equation}	\label{eq1016wed}
rE:=\{(r^2t,rx):(t,x)\in E\}.
\end{equation}
Then
\begin{equation}	\label{eq1717mon}
\capacity(rE)=r^d\capacity(E).
\end{equation}
Every singleton has zero thermal capacity.
\end{lemma}

\begin{lemma}	\label{lem02}
Let $(E_j)$ be a sequence of compact subsets of $\mathbb R^{d+1}$, and set $E=\bigcup_{j=1}^\infty E_j$.
Then
\[
\red_E\le\sum_{j=1}^\infty\red_{E_j} \qquad\text{in }\,\mathbb R^{d+1}.
\]
Moreover, if $F$ is compact and $(K_m)$ is a decreasing sequence of compact neighborhoods of $F$ with $\bigcap_m K_m=F$, then
\[
\red_{K_m}\downarrow\red_F \qquad\text{in }\,\mathbb R^{d+1}.
\]
\end{lemma}

\subsection{Boundary regularity via reduced functions}
For a compact set $K$, define its \emph{regularized reduced function} (or \emph{balayage}) by
\[
\bal_K(z) =\lim_{r\downarrow0}\inf_{Q_r(z)}\red_K.
\]
By construction, $\bal_K$ is lower semicontinuous.

\begin{lemma}\label{lem05}
Let $K\subset\mathbb R^{d+1}$ be compact. Then
\[
0\le\bal_K\le\red_K\le1, \qquad \red_K=1\quad\text{on }K, \qquad \bal_K=1\quad\text{on }K^\circ.
\]
Moreover, $\bal_K$ is $P$-superparabolic in $\mathbb R^{d+1}$, and $\red_K=\bal_K$ is $P$-caloric in $K^\complement$.
\end{lemma}

\begin{proposition}	\label{prop01}
Let $\Omega\subset\mathbb R^{d+1}$ be bounded and open, and let $z_0\in\partial\Omega$.
For $r>0$, set
\[
K_r=\Omega^\complement\cap\overline Q_r(z_0).
\]
Then $z_0$ is regular if and only if
\[
\bal_{K_r}(z_0)=1 \quad\text{for all }r>0.
\]
\end{proposition}
\begin{proof}
Suppose first that $z_0$ is regular.
Fix $r>0$ and choose $f\in C(\partial\Omega)$ such that $0\le f\le1$, $f(z_0)=1$, and $\supp f \subset Q_r(z_0)$.
If $v$ is nonnegative and $P$-superparabolic in $\mathbb R^{d+1}$ with $v\ge 1$ on $K_r$, then $v\vert_\Omega\in\mathcal U_f$, the upper Perron class.
Consequently,
\[
\overline H_f \le \red_{K_r}\qquad\text{in }\Omega.
\]
Since $z_0$ is regular and $0\le \red_{K_r}\le 1$,
\[
\lim_{\Omega\ni z\to z_0}\red_{K_r}(z)=1.
\]
Moreover, $\red_{K_r}=1$ on $\Omega^\complement\cap Q_r(z_0)$.
It follows that $\bal_{K_r}(z_0)=1$.

Conversely, suppose that $\bal_{K_r}(z_0)=1$ for all $r>0$.
We first show that there exists a constant $M>2$ such that, for every $0<r<1/M$, one can find a nonnegative $P$-superparabolic function $w_r$ satisfying
\begin{equation}	\label{eq1503sun}
w_r\ge1\quad\text{on }\overline Q_r(z_0), \qquad
w_r\le\tfrac14\quad\text{on }\partial Q_{Mr}(z_0).
\end{equation}

In this construction, we use the Gaussian bounds for $0<t-\tau<2$.
Choose $\eta\in C_c^\infty((-3,-2))$ such that $\eta\ge0$ and $\int\eta=1$, and $\phi\in C_c^\infty(B_2)$ such that $0\le\phi\le 1$ and $\phi=1$ on $B_1$.
Writing $z_0=(t_0,x_0)$, define
\[
g_r(t,x)=r^{-2}\eta\left(\frac{t-t_0}{r^2}\right) \phi\left(\frac{x-x_0}{r}\right),
\qquad w_r=A\mathscr Gg_r,
\]
where $A>0$ will be chosen below.
Since $g_r\in C_c^\infty(\mathbb R^{d+1})$, standard parabolic regularity implies that $w_r$ is continuous.
Moreover, $w_r\ge0$ and $Pw_r=Ag_r\ge0$, so $w_r$ is $P$-superparabolic.

For $(t,x)\in\overline Q_r(z_0)$ and $\tau\in\supp\eta((\cdot-t_0)/r^2)$, we have
\[
r^2\le t-\tau\le4r^2,\qquad \abs{x-x_0} \le r.
\]
Since $\abs{x-y} \le 2r$ for $y\in B_r(x_0)$, the Gaussian lower bound gives
\[
\Gamma(t,x,\tau,y) \ge c_0 2^{-d}e^{-4b_0}r^{-d}.
\]
Consequently,
\[
\begin{aligned}
\int_{\mathbb R^d}\Gamma(t,x,\tau,y) \phi\left(\frac{y-x_0}{r}\right)\,dy &\ge\int_{B_r(x_0)}\Gamma(t,x,\tau,y)\,dy\\
&\ge c_0 2^{-d}e^{-4b_0}\abs{B_1}=:\kappa>0.
\end{aligned}
\]
Since $\eta\ge 0$ and $\int\eta=1$, we obtain $w_r\ge A\kappa$ on $\overline Q_r(z_0)$. Thus, we fix $A=\kappa^{-1}$.

Next, let $M>2$.
When $t=t_0+M^2r^2$, we have $t-\tau\ge(M^2+2)r^2$ for $\tau$ in the temporal support of $g_r$.
Hence, on this part of $\partial Q_{Mr}(z_0)$, the Gaussian upper bound gives
\[
w_r(t,x)\le AC_0 2^d \abs{B_1}(M^2+2)^{-d/2}.
\]
On the lateral boundary, where $\abs{x-x_0}=Mr$, we have $\abs{x-y}\ge(M-2)r$ for $y$ in the spatial support of $g_r$.
For $\tau<t$, setting $\theta=(t-\tau)/r^2$, we obtain
\[
\begin{aligned}
w_r(t,x) &\le AC_0 2^d\abs{B_1}\sup_{\theta>0} \theta^{-d/2}\exp\left(-a_0\frac{(M-2)^2}{\theta}\right)\\
&=AC_0 2^d\abs{B_1} \left(\frac{d}{2ea_0}\right)^{d/2}(M-2)^{-d}.
\end{aligned}
\]
When $t=t_0-M^2r^2$, we have $w_r=0$ since $M^2>3$.
Choose $M$ sufficiently large that both upper bounds are at most $1/4$.
This proves \eqref{eq1503sun}.

Since $w_r\ge1$ on $K_r$, the definition of $\red_{K_r}$ gives
\[
\bal_{K_r}\le\red_{K_r}\le w_r.
\]
Set
\[
v_r=1-\bal_{K_r}.
\]
By Lemma~\ref{lem05}, $v_r$ is nonnegative and $P$-caloric in $\Omega$.
Moreover, since $\bal_{K_r}(z_0)=1$ and $\bal_{K_r}$ is lower semicontinuous,
\[
\lim_{\Omega\ni z\to z_0}v_r(z)=0.
\]
On the other hand, \eqref{eq1503sun} and the continuity of $w_r$ imply that $v_r \ge 1/2$ in a neighborhood of $\Omega\cap\partial Q_{Mr}(z_0)$.

Let $f\in C(\partial\Omega)$.
Subtracting $f(z_0)$ and scaling, we may assume $f(z_0)=0$ and $\abs{f} \le 1$.
Given $\epsilon\in(0,1)$, choose $r>0$ sufficiently small that
\[
\abs{f} \le\epsilon \quad\text{on }\partial\Omega\cap\overline Q_{Mr}(z_0).
\]
Define
\[
v(z)=
\begin{cases}
\min\{1,\epsilon+2v_r(z)\},
&z\in\Omega\cap Q_{Mr}(z_0),\\
1,&z\in\Omega\setminus Q_{Mr}(z_0).
\end{cases}
\]
Since $v_r\ge 1/2$ near $\Omega\cap\partial Q_{Mr}(z_0)$, the two expressions agree there.
Thus $v$ is $P$-superparabolic in $\Omega$.
Moreover, $v\in\mathcal U_f$, and hence
\[
\overline H_f \le v.
\]
Since $v_r(z)\to0$ as $z\to z_0$ within $\Omega$,
\[
\limsup_{\Omega\ni z\to z_0}\overline H_f(z)\le\epsilon.
\]
Applying the same argument to $-f$, we obtain
\[
\liminf_{\Omega\ni z\to z_0}\underline H_f(z)\ge-\epsilon.
\]
Since $\underline H_f\le\overline H_f$, letting $\epsilon\downarrow 0$ gives
\[
\lim_{\Omega\ni z\to z_0}\overline H_f(z)=0.
\]
Therefore $z_0$ is regular.
\end{proof}

\begin{lemma}	\label{lem06}
Let $K\subset\mathbb R^{d+1}$ be compact and $z_0=(t_0,x_0)$.
For $n\ge 1$, set
\[
K_n=K\cap\{t\le t_0-1/n\}.
\]
Then
\[
\bal_K(z_0) =\sup_{n\ge1}\red_{K_n}(z_0) =\red_{K\cap\{t<t_0\}}(z_0).
\]
\end{lemma}

\begin{proof}
We first observe that every bounded nonnegative $P$-superparabolic function $v$ satisfies
\begin{equation}	\label{eq1418sun}
v(t_0,x_0) =\liminf_{\substack{(t,x)\to(t_0,x_0)\\t<t_0}}v(t,x).
\end{equation}

By lower semicontinuity,
\[
v(t_0,x_0)\le \liminf_{\substack{(t,x)\to(t_0,x_0)\\t<t_0}}v(t,x).
\]
If the inequality were strict, there would exist $a>v(t_0,x_0)$ and $r>0$ such that
\[
v\ge a\quad\text{in } (t_0-r^2,t_0)\times B_r(x_0).
\]
Choose $\phi\in C_c^\infty(B_r(x_0))$ with $0 \le \phi \le 1$ and $\phi(x_0)=1$.
For $s<t_0$ sufficiently close to $t_0$, let $h_s$ be the $P$-caloric function in $(s,t_0+r^2)\times B_r(x_0)$ with initial datum $a\phi$ and zero lateral boundary data.
Comparison gives $v\ge h_s$.
Since $h_s(t_0,x_0)\to a$ as $s\uparrow t_0$, we obtain $v(t_0,x_0)\ge a$, a contradiction.
This proves \eqref{eq1418sun}.

For $\tau\in\mathbb R$, the function $\mathbbm1_{\{t>\tau\}}$ is nonnegative and $P$-superparabolic.
If $v$ is $P$-superparabolic with $0\le v\le 1$ and $v=1$ on $K\cap\{t\le\tau\}$, then $v+\mathbbm1_{\{t>\tau\}}$ is nonnegative and $P$-superparabolic, with $v+\mathbbm1_{\{t>\tau\}}\ge 1$ on $K$.
Taking the infimum over such $v$, we obtain
\[
\red_K\le \red_{K\cap\{t\le\tau\}} \qquad\text{on }\,\{t\le\tau\}.
\]
The reverse inequality follows from monotonicity.
Therefore,
\[
\red_K=\red_{K\cap\{t\le\tau\}} \qquad\text{on }\,\{t\le\tau\}.
\]

For every $(t,x)$ with $t<t_0$, we may choose $n$ sufficiently large that $t<t_0-1/n$.
The preceding equality then holds in a neighborhood of $(t,x)$, and hence
\[
\bal_K(t,x) =\bal_{K_n}(t,x).
\]

Set
\[
w=\sup_{n\ge1}\bal_{K_n}.
\]
Since $K_n$ increases with $n$, Lemmas \ref{lem03} and \ref{lem05} imply that $w$ is bounded and $P$-superparabolic.
Moreover,
\[
w=\bal_K \qquad\text{on }\,\{t<t_0\}.
\]
Applying \eqref{eq1418sun} to both functions at $z_0$, we conclude that
\[
\bal_K(z_0)=w(z_0) =\sup_{n\ge1}\bal_{K_n}(z_0).
\]
Since $z_0\notin K_n$, Lemma~\ref{lem05} gives
\[
\bal_{K_n}(z_0) =\red_{K_n}(z_0).
\]

Finally, every compact subset of $K\cap\{t<t_0\}$ is contained in $K_n$ for sufficiently large $n$.
Therefore, by the definition of $\red_E$ for arbitrary sets,
\[
\sup_{n\ge 1}\red_{K_n}(z_0) =\red_{K\cap\{t<t_0\}}(z_0). \qedhere
\]
\end{proof}

\begin{proposition}	\label{prop02}
Let $\Omega\subset\mathbb R^{d+1}$ be bounded and open, and let $z_0\in\partial\Omega$. Set
\[
E_\epsilon=E_\epsilon(z_0):=\Omega^\complement\cap Q_\epsilon^-(z_0).
\]
Then $z_0$ is regular if and only if
\[
\red_{E_\epsilon}(z_0)=1\qquad\text{for every }\epsilon>0.
\]
In particular, if $z_0$ is irregular, then
\[
\red_{E_\epsilon}(z_0)<1
\]
for all sufficiently small $\epsilon>0$.
\end{proposition}

\begin{proof}
Let $K_r=\Omega^\complement\cap\overline Q_r(z_0)$.
By Proposition~\ref{prop01} and Lemma~\ref{lem06}, $z_0=(t_0,x_0)$ is regular if and only if
\[
\red_{K_r\cap\{t<t_0\}}(z_0)=1
\]
for every $r>0$. Since
\[
K_{\epsilon/2}\cap\{t<t_0\} \subset E_\epsilon \subset K_\epsilon\cap\{t<t_0\},
\]
the equivalence follows from monotonicity of reduced functions.

If $z_0$ is irregular, then $\red_{K_r\cap\{t<t_0\}}(z_0)<1$ for some $r>0$.
Since $E_\epsilon\subset K_r\cap\{t<t_0\}$ whenever $0<\epsilon\le r$, monotonicity gives the last assertion.
\end{proof}

By Proposition~\ref{prop02}, boundary regularity is a local property and is thus defined in the same way for unbounded open sets by localization.
\section{Auxiliary estimates}	\label{sec3}
\subsection{A heat-ball estimate at the critical exponent}
The following lemma is an analogue of the strong Harnack inequality of Evans and Gariepy \cite[Lemma~3.2]{EG1982}, which plays a key role in their proof of the Wiener criterion for the heat equation.
Our estimate holds under the two-sided Gaussian bounds \eqref{eq0157sat}, with heat balls defined using the critical exponent $b_0$ in the Gaussian lower bound.
Its proof relies on Gaussian estimates and comparison principles rather than the explicit structure of the heat kernel.

For $q>0$, define
\[
\cE_q=\{(t,x)\in\mathbb R^{d+1}:t<0,\ G_{b_0}(-t,x)>q\}.
\]
\begin{lemma}	\label{lem10}
Let $\alpha\ge 1$.
There exist constants $\beta\in(0,1)$ and $\gamma>0$, depending only on $d$, $\alpha$, and the constants in Condition~\ref{gaussian}, such that the following holds.
For every $q \ge 1$ with $r=q^{-1/d}$,
\[
\{(-8r^2,y):\abs{y}\le\alpha r\}\subset \cE_{\beta q}.
\]
Moreover, every bounded nonnegative $P$-caloric function $u$ in $\cE_{\beta q}$ satisfies
\[
\inf_{\cE_q}u \ge\gamma r^{-d}\int_{\abs{y}<\alpha r}u(-8r^2,y)\,dy.
\]
\end{lemma}

\begin{proof}
By the observation following Condition~\ref{gaussian}, we may assume $T_0>16$.
Since $r\le1$, parabolic scaling \eqref{eq1016wed} reduces the proof to the case $q=1$.

Let $M>\frac{d}{2}\log 8$ be a parameter to be fixed later, and set $\beta=e^{-M}$.
Write $\tau=-t$, and define
\[
R_M(\tau) =\sqrt{\frac{\tau}{b_0} \left(M+\frac{d}{2}\log\frac{1}{\tau}\right)} \qquad (0<\tau \le 8).
\]
Choose
\[
M>\max\left\{\frac{d}{2}(1+\log8),\,\frac{d}{2}\log7+\frac{b_0}{7}(\alpha+1)^2\right\},
\]
to be increased if necessary below.
Then $R_M$ is increasing and concave on $(0,8]$, and $R_M(7)>\alpha+1$.
In particular, $R_M(8)>\alpha$, so
\[
\{(-8,y):\abs{y}\le\alpha\}\subset \cE_\beta.
\]

For $0\le s<8$, define
\[
\cD_s:=\{(-\tau,x):s<\tau<8,\ \abs{x}<R_M(\tau)\}.
\]
In particular, $\cD_0=\cE_\beta\cap\{-8<t<0\}$.
The lateral boundary of $\cD_s$ is
\[
\partial_x \cD_s :=\{(-\tau,x):s<\tau<8,\ \abs{x}=R_M(\tau)\}.
\]

Observe that
\[
\cE_1=\{(-s,x):0<s<1,\ |x|<R_0(s)\},
\qquad
R_0(s):=\sqrt{\frac{ds}{2b_0}\log\frac1s}.
\]
Thus, for $z^\circ=(-s,x^\circ)\in \cE_1$, we have
\begin{equation}	\label{eq1029fri}
\abs{x^\circ}<R_0(s)\le\sqrt{\frac{d}{2eb_0}}.
\end{equation}

We first construct a nonnegative $P$-superparabolic function $v$ in $\cD_s$ such that
\begin{equation}	\label{eq1030fri}
\liminf_{\cD_s \ni z\to\zeta} v(z) \ge1 \quad (\zeta\in\partial_x \cD_s),\qquad \limsup_{\cD_s \ni z\to z^\circ} v(z) \le Ce^{-cM},
\end{equation}
where $C,c>0$ depend only on $d$, $a_0$, $b_0$, and $C_0$.

For $h,a>0$ with $h\le 8$, choose $\psi\in C^\infty(\mathbb R^d)$ such that
\[
0 \le \psi \le 1,\qquad \psi(x)=0\;\text{ for }\;\abs{x-x^\circ}<a/2,\qquad \psi(x)=1\;\text{ for }\; \abs{x-x^\circ} \ge 3a/4.
\]
Put $\eta=-s-2h$ and
\[
w_{h,a}(t,x)=
\begin{cases}
2\mathscr S_\eta\psi(t,x),&t>\eta,\\
0,&t\le\eta.
\end{cases}
\]
This function is nonnegative and $P$-superparabolic in $\mathbb{R}^{d+1}$.
For $-s-h\le t\le-s$ and $\abs{x-x^\circ} \ge a$, the mass conservation identity and the Gaussian upper bound give
\[
w_{h,a}(t,x) =2-2\mathscr S_\eta(1-\psi)(t,x)\ge 2(1-Ce^{-ca^2/h}).
\]
Similarly, since $\psi=0$ on $B_{a/2}(x^\circ)$,
\[
w_{h,a}(z^\circ)\le Ce^{-ca^2/h}.
\]
Thus, if $a^2/h$ is sufficiently large, then
\begin{equation}	\label{eq1032fri}
\begin{aligned}
w_{h,a}(t,x)&\ge 1 &&\text{if }\;-s-h\le t\le-s,\; \abs{x-x^\circ} \ge a,\\
w_{h,a}(z^\circ)&\le Ce^{-ca^2/h}.
\end{aligned}
\end{equation}

Writing $\tau=-t$, we divide $\partial_x \cD_s$ into two parts.
First consider $s<\tau\le2s$. Set
\[
\delta_s=R_M(s)-R_0(s),\qquad \Delta_s=R_M(2s)-R_M(s),
\]
and write
\[
A_s=M+\frac{d}{2}\log\frac{1}{s}.
\]
Since $0<s<1$, we have $A_s> M$, and
\[
\delta_s =\sqrt{\frac{s}{b_0}} \left(\sqrt{A_s}-\sqrt{A_s-M}\right) \ge \sqrt{\frac{s}{b_0}}\,\frac{M}{2\sqrt{A_s}}.
\]
On the other hand, if $M \ge 2d\log2$, then $d\log2\le A_s/2$, and hence
\[
\Delta_s =\sqrt{\frac{s}{b_0}} \left(\sqrt{2A_s-d\log2}-\sqrt{A_s}\right)
\ge \left(\sqrt{\frac32}-1\right) \sqrt{\frac{s}{b_0}}\sqrt{A_s}.
\]
Consequently,
\[
\frac{\delta_s\Delta_s}{s} \ge c M.
\]

By concavity of $R_M$,
\[
R_M(\tau)\ge R_M(s) +\frac{\tau-s}{s}\,\Delta_s, \qquad s\le\tau\le 2s.
\]
Thus, for $(-\tau,x)\in\partial_x \cD_s$ with $s<\tau\le2s$, the preceding inequality and \eqref{eq1029fri} yield
\[
\abs{x-x^\circ} \ge R_M(\tau)-\abs{x^\circ} \ge\delta_s+\frac{\tau-s}{s}\Delta_s.
\]

For $n\ge0$, set
\[
h_n=2^{-n}s,\qquad a_n=\delta_s+2^{-n-1}\Delta_s.
\]
By the preceding estimate, if $(-\tau,x)\in\partial_x \cD_s$ with $s+h_n/2 \le \tau \le s+h_n$, then $\abs{x-x^\circ} \ge a_n$.
Moreover,
\[
\frac{a_n^2}{h_n} =\frac{\delta_s^2}{s}\,2^n +\frac{\delta_s\Delta_s}{s}+\frac{\Delta_s^2}{4s}\,2^{-n} \ge cM.
\]
Thus, taking $M$ sufficiently large, we obtain $w_{h_n,a_n}\ge 1$ on this portion of $\partial_x\cD_s$ by \eqref{eq1032fri}.
To estimate their sum at $z^\circ$, set
\[
\theta_s=\frac{\delta_s\Delta_s}{2s}, \qquad \nu_s=\log_2\frac{\Delta_s}{2\delta_s}.
\]
Then
\[
\frac{a_n^2}{h_n} =\theta_s\left(2^{n-\nu_s}+2+2^{-(n-\nu_s)}\right).
\]
For each integer $m\ge0$, there are at most two integers $n\ge0$ satisfying
\[
m\le \abs{n-\nu_s}<m+1.
\]
For such $n$,
\[
2^{n-\nu_s}+2^{-(n-\nu_s)}\ge2^m.
\]
Since $\theta_s\ge cM$, it follows that
\[
\sum_{n=0}^\infty w_{h_n,a_n}(z^\circ) \le C\sum_{n=0}^\infty e^{-ca_n^2/h_n} \le 2Ce^{-c\theta_s} \sum_{m=0}^\infty e^{-c\theta_s2^m} \le Ce^{-cM}.
\]
Set
\[
v_1=\sum_{n=0}^\infty w_{h_n,a_n}.
\]
The sum is locally finite in $\cD_s$.
Moreover, since $\delta_s>0$, the estimate
\[
w_{h_n,a_n}(t,x)\le Ce^{-c\delta_s^2/h_n}
\]
for $t\le-s$ and $\abs{x-x^\circ} \le\delta_s/8$ shows that the sum converges uniformly near $z^\circ$ from within $\cD_s$.
Consequently, for every $(-\tau,x)\in\partial_x \cD_s$ with $s<\tau\le2s$, we have
\[
\liminf_{\cD_s\ni z\to (-\tau,x)}v_1(z)\ge 1.
\]
Moreover,
\[
\limsup_{\cD_s\ni z\to z^\circ}v_1(z) \le Ce^{-cM}.
\]

Next consider $2s<\tau<8$. Since $M\ge\frac d2\log2$,
we have
\[
R_M(2s)^2\ge 2R_0(s)^2.
\]
Thus, by \eqref{eq1029fri} and the monotonicity of $R_M$,
\[
\abs{x^\circ} \le2^{-1/2}R_M(\tau), \qquad 2s<\tau<8.
\]
For integers $m\ge 1$ with $2^ms<8$, set
\[
h_m=\min\{2^{m+1}s,8\}-s,\qquad
a_m=(1-2^{-1/2})R_M(2^ms).
\]
If $(-\tau,x)\in\partial_x \cD_s$ with $2^ms<\tau\le\min\{2^{m+1}s,8\}$, then $\abs{x-x^\circ} \ge a_m$.
Moreover, since $h_m\le2^{m+1}s$,
\[
\frac{a_m^2}{h_m} \ge c\left(M+\frac d2\log\frac1{2^ms}\right).
\]
Since $2^ms<8$, the preceding lower bound for $a_m^2/h_m$ is sufficiently large when $M$ is large.
Thus, \eqref{eq1032fri} gives $w_{h_m,a_m}\ge 1$ on the portion of $\partial_x \cD_s$ with $2^ms<\tau\le\min\{2^{m+1}s,8\}$ and
\[
w_{h_m,a_m}(z^\circ) \le Ce^{-ca_m^2/h_m} \le Ce^{-cM}(2^ms)^\sigma
\]
for some $\sigma>0$.
Set
\[
v_2=\sum_{\substack{m\ge 1\\2^ms<8}} w_{h_m,a_m}.
\]
Then, for every $(-\tau,x)\in\partial_x \cD_s$ with $2s<\tau<8$, we have
\[
\liminf_{\cD_s\ni z\to (-\tau, x)}v_2(z)\ge 1.
\]
Moreover,
\[
\limsup_{\cD_s\ni z\to z^\circ}v_2(z) \le Ce^{-cM} \sum_{\substack{m\ge 1\\2^ms<8}} (2^ms)^\sigma
\le Ce^{-cM} \frac{8^\sigma}{1-2^{-\sigma}} \le Ce^{-cM}.
\]
Thus $v=v_1+v_2$ is a nonnegative $P$-superparabolic function in $\cD_s$ satisfying \eqref{eq1030fri}.

We next compare $u$ with the function
\[
\tilde u(t,x)=\int_{\abs{y}<\alpha} \Gamma(t,x,-8,y)u(-8,y)\,dy.
\]
Set
\[
\mu:=\int_{\abs{y}<\alpha}u(-8,y)\,dy.
\]
Since $7<8-s<8$, by \eqref{eq1029fri} and the Gaussian lower bound,
\[
\tilde u(z^\circ)\ge c_1 \mu,\qquad\text{where }\;
c_1=c_0\,8^{-d/2} \exp\left\{-\frac{b_0}{7}\left(\sqrt{\frac{d}{2eb_0}}+\alpha\right)^2\right\}.
\]

On the lateral boundary of $\cD_s$, we have
\begin{equation}	\label{eq1036fri}
\tilde u\le c_2 \mu,\qquad\text{where }\;
c_2=C_0\max\left\{1,\left(\frac{d}{2ea_0}\right)^{d/2}\right\}.
\end{equation}
Indeed, by the Gaussian upper bound,
\[
\Gamma(t,x,-8,y)\le C_0G_{a_0}(t+8,x-y).
\]
When $t\ge-7$, the Gaussian kernel is bounded by $1$.
When $-8<t<-7$, we have $\abs{x}=R_M(-t)>\alpha+1$, and hence $\abs{x-y}\ge1$ for $\abs{y}<\alpha$.
The desired bound follows from
\[
\sup_{0<\theta<1}\theta^{-d/2}e^{-a_0/\theta} \le\left(\frac{d}{2ea_0}\right)^{d/2}.
\]
The initial-time comparison $\tilde u\le u$ at $t=-8$ is justified by approximating
$\mathbbm1_{\{\abs{y}<\alpha\}}$ from below by smooth cutoffs.

By \eqref{eq1030fri}, \eqref{eq1036fri}, and the initial-time comparison, the function
$u-\tilde u+c_2\mu v$ is $P$-superparabolic in $\cD_s$ and satisfies
\[
\liminf_{\cD_s\ni z\to\zeta} (u(z)-\tilde u(z)+c_2\mu v(z))\ge 0
\]
for every $\zeta\in\partial \cD_s\cap\{t<-s\}$.

The comparison principle therefore yields $u\ge\tilde u-c_2\mu v$ in $\cD_s$.
Letting $z\to z^\circ$ within $\cD_s$, and recalling the definition of $\mu$, we obtain
\[
u(z^\circ) \ge \tilde u(z^\circ) -c_2\mu\limsup_{\cD_s\ni z\to z^\circ}v(z)
\ge \left(c_1-c_2Ce^{-cM}\right)\int_{\abs{y}<\alpha}u(-8,y)\,dy.
\]
Choose $M$ sufficiently large that $c_2Ce^{-cM}\le c_1/2$.
Then
\[
u(z^\circ)\ge\frac{c_1}{2} \int_{\abs{y}<\alpha}u(-8,y)\,dy.
\]
Since $z^\circ \in \cE_1$ was arbitrary, this proves the assertion for $q=1$.

For general $q\ge1$, set $r=q^{-1/d}$ and apply the normalized estimate to $u(r^2t,rx)$.
The rescaled fundamental solution satisfies the same Gaussian bounds for time differences up to $T_0/r^2\ge16$.
The conclusion follows with $\gamma=c_1/2$ and $\beta=e^{-M}$.
\end{proof}

\subsection{Penalization and a capacitary estimate}
The following approximation is closely related to the penalization method for obstacle problems.
\begin{lemma}	\label{lem11}
Let $K$ be a compact subset of $(-\infty,0)\times\mathbb R^d$ and let $\chi\in C_c^\infty(K^\circ)$ satisfy $0\le\chi\le 1$.
Then there exist bounded nonnegative functions $f_n$ with $\supp f_n\subset\supp\chi$ such that $u_n=\mathscr Gf_n$ satisfies
\[
0\le u_n\le 1,\qquad Pu_n=f_n,\qquad u_n \le\red_K.
\]
Moreover, if $\chi=1$ in a neighborhood of a compact set $E\subset K^\circ$, then $u_n\to1$ uniformly on $E$.
\end{lemma}
\begin{proof}
Choose $s<0$ such that $K\subset(s,0)\times\mathbb R^d$.
Let $X$ be the space of bounded measurable functions on $[s,0]\times\mathbb R^d$ taking values in $[0,1]$, equipped with the supremum norm.
For each $n\in\mathbb N$, define $\mathcal T_n$ on $X$ by
\[
(\mathcal T_n v)(t,x) =e^{-n(t-s)} +n\int_s^t e^{-n(t-\tau)} \int_{\mathbb R^d} \Gamma(t,x,\tau,\xi)(1-\chi(\tau,\xi)) v(\tau,\xi)\,d\xi\,d\tau.
\]
By positivity and mass conservation of $\Gamma$, $\mathcal T_n$ maps $X$ into itself and
\[
\norm{\mathcal T_n v-\mathcal T_n w}_\infty \le(1-e^{ns})\,\norm{v-w}_\infty.
\]
Thus, by the Banach fixed-point theorem, there exists a unique $v_n\in X$ such that $v_n=\mathcal T_n v_n$.
By Duhamel's formula, $v_n$ satisfies
\[
Pv_n+n v_n=n(1-\chi)v_n,\qquad v_n(s,\cdot)=1,
\]
or equivalently,
\[
Pv_n+n\chi v_n=0,\qquad v_n(s,\cdot)=1.
\]

Define $f_n=n\chi v_n$ on $[s,0]\times\mathbb R^d$ and extend $f_n=0$ elsewhere. Set
\[
u_n=\mathscr G f_n.
\]
Since $P(1-v_n)=f_n$ and $(1-v_n)(s,\cdot)=0$, uniqueness of bounded solutions to the Cauchy problem gives
\[
u_n=1-v_n\qquad\text{on }\;[s,0]\times\mathbb R^d.
\]
Hence $0\le u_n\le1$ on $[s,0]\times\mathbb R^d$, and $u_n=0$ on $(-\infty,s]\times\mathbb R^d$.
For $t>0$, the reproduction identity gives
\[
u_n(t,x)=\int_{\mathbb R^d} \Gamma(t,x,0,\xi)u_n(0,\xi)\,d\xi.
\]
By positivity and mass conservation of $\Gamma$, we also have $0\le u_n\le1$ for $t>0$.
Consequently,
\[
0\le u_n\le1,\qquad Pu_n=f_n,\qquad \supp f_n\subset\supp\chi.
\]

We next show that $u_n \le \red_K$.
Let $v$ be any nonnegative $P$-superparabolic function in $\mathbb R^{d+1}$ satisfying $v\ge 1$ on $K$.
Fix an arbitrary $t_1>0$.
For sufficiently large $R$, set
\[
D_R=((s,t_1)\times B_R)\setminus K.
\]
Since $f_n$ is supported in $K$, the function $v-u_n$ is $P$-superparabolic in $D_R$.

Since $f_n$ is bounded and compactly supported, the Gaussian upper bound for $0<t-\tau<t_1-s$ implies that
\[
\epsilon_R:=\sup_{(s,t_1)\times\partial B_R}u_n \to 0\qquad(R\to\infty).
\]
Since $u_n=0$ on $\{s\}\times\overline B_R$, $u_n\le1\le v$ on $K$, and $u_n\le\epsilon_R$ on $(s,t_1)\times\partial B_R$, we have
\[
\liminf_{D_R\ni z\to\zeta}\, (v(z)+\epsilon_R-u_n(z))\ge0 \qquad(\zeta\in\partial D_R\cap\{t<t_1\}).
\]
Here we use the lower semicontinuity of $v$ and the continuity of $u_n$.
The comparison principle then gives
\[
u_n\le v+\epsilon_R\qquad\text{in }\,D_R.
\]
Since $u_n\le v$ on $K$, the same inequality holds in $(s,t_1)\times B_R$.
Letting $R\to\infty$ and taking the infimum over all such $v$, we obtain $u_n\le\red_K$ in $(s,t_1)\times\mathbb R^d$.
Since $t_1>0$ is arbitrary and $u_n=0$ for $t\le s$, the inequality holds in $\mathbb R^{d+1}$.

Finally, suppose that $\chi=1$ in a neighborhood of a compact set $E\subset K^\circ$.
Write $z^\circ=(t^\circ,x^\circ)$ for a point of $E$.
By compactness, there exist $r,h_0>0$ such that
\[
[t^\circ-h_0,t^\circ]\times\overline B_r(x^\circ) \subset\{\chi=1\} \qquad(z^\circ\in E),
\]
and $t^\circ-h_0>s$ for all $z^\circ\in E$.
With $r$ fixed, choose $0<h\le h_0$ sufficiently small that $r^2/h$ exceeds the threshold required in \eqref{eq1032fri}.

For each $z^\circ=(t^\circ,x^\circ)\in E$, set
\[
Q=(t^\circ-h,t^\circ)\times B_r(x^\circ),
\]
and let $w$ be the function constructed in \eqref{eq1032fri}, centered at $z^\circ$ with parameters $h$ and $r$.
Then $w\ge 0$, $Pw=0$ in $Q$, and
\[
w\ge 1\quad\text{on }(t^\circ-h,t^\circ) \times\partial B_r(x^\circ), \qquad w(z^\circ)\le Ce^{-cr^2/h}.
\]

Define
\[
\tilde w(t,x)=e^{-n(t-t^\circ+h)}+w(t,x).
\]
Since $\chi=1$ and $Pw=0$ in $Q$, we have
\[
P\tilde w+n\tilde w=nw\ge0 \qquad\text{in }Q.
\]
Moreover, $\tilde w\ge1\ge v_n$ on $\partial_p Q$.
Since $Pv_n+nv_n=0$ in $Q$, comparison and passage to $z^\circ$ from within $Q$ give
\[
v_n(z^\circ)\le e^{-nh}+Ce^{-cr^2/h}.
\]
The estimate is uniform for $z^\circ\in E$. Hence
\[
\limsup_{n\to\infty}\sup_E v_n
\le Ce^{-cr^2/h}.
\]
Letting $h\downarrow0$ gives $\sup_E v_n\to0$, and therefore $u_n\to 1$ uniformly on $E$.
\end{proof}

\begin{lemma}		\label{lem14}
Let $E$ be a nonempty compact subset of $(-1,0)\times\mathbb R^d$. Then
\[
\red_E(0) \le C \left( \max_{(t,x)\in E}G_{a_0}(-t,x)\right)\,\capacity(E),
\]
where $C>0$ depends only on $d$ and the constants in Condition~\ref{gaussian}.
\end{lemma}
\begin{proof}
Choose decreasing compact neighborhoods $K_m\Subset(-1,0)\times\mathbb R^d$ of $E$ such that $K_m\downarrow E$.
For each $m$, choose $\chi\in C_c^\infty(K_m^\circ)$ with $0\le\chi\le1$ and $\chi=1$ near $E$.
Let $f_n$ and $u_n=\mathscr Gf_n$ be as in Lemma~\ref{lem11}.

Set $w_n=G_{b_0} f_n$.
Thus, the Gaussian lower bound and $u_n\le 1$ give
\[
c_0w_n(t,x)\le u_n(t,x)\le1 \qquad(t\le0).
\]
For $t>0$, the reproduction identity for the heat kernel gives
\[
w_n(t,x)= \left(\frac{b_0}{\pi}\right)^{d/2} \int_{\mathbb R^d} G_{b_0}(t,x-y)w_n(0,y)\,dy.
\]
Since $w_n(0,\cdot)\le c_0^{-1}$, positivity and mass conservation yield
\[
w_n\le c_0^{-1}\qquad\text{in }\mathbb R^{d+1}.
\]

Since $f_n\ge 0$ is supported in $K_m$ and $G_{b_0}f_n\le c_0^{-1}$ in $\mathbb R^{d+1}$, the definition of Gaussian capacity gives
\[
\int_{\mathbb R^{d+1}}f_n \le c_0^{-1}\mathcal C_{b_0}(K_m).
\]
Set
\[
M_m:=\max_{(t,x)\in K_m}G_{a_0}(-t,x).
\]
The Gaussian upper bound and Lemma~\ref{lem01} then give
\[
u_n(0)\le C_0M_m \int_{\mathbb R^{d+1}} f_n \le\frac{C_0}{c_0}\,M_m \mathcal C_{b_0}(K_m) \le C M_m \capacity(K_m).
\]

Since $u_n\to 1$ uniformly on $E$, we have
\[
\delta_n:=1-\min_E u_n\to 0.
\]
For all sufficiently large $n$, the function $u_n/(1-\delta_n)$ is nonnegative and $P$-superparabolic, with $u_n/(1-\delta_n)\ge 1$ on $E$.
Thus, by the definition of the reduced function, $\red_E \le u_n/(1-\delta_n)$.
Consequently,
\[
\red_E(0) \le\liminf_{n\to\infty}u_n(0) \le CM_m \capacity(K_m).
\]
Finally, as $m\to\infty$, we have
\[
M_m\to\max_{(t,x)\in E}G_{a_0}(-t,x),
\qquad
\capacity(K_m)\to\capacity(E)
\]
by continuity and Lemma~\ref{lem01}, respectively.
\end{proof}

\section{Qualitative Wiener-type criteria and applications}\label{sec4}
For divergence form operators, Condition~\ref{gaussian} with $T_0=\infty$ follows from Aronson's estimates \cite{Aronson}.
For non-divergence form operators with $\mathbf A\in\mathrm{DMO}_x$, it follows from \cite{DEK21}.

Fix $0<\rho<1$.
By translation, we may assume $z_0=(0,0)=0$.
For simplicity, write
\[
\Omega_k^\ell=\Omega_k^\ell(0,\rho).
\]

\subsection{Sufficiency via Gaussian level sets}
\begin{lemma}	\label{lem12}
For any $0<\delta<1$, there exists a constant $C>0$ such that the following holds.
Let $K$ be a compact subset of $(-1,0)\times\mathbb R^d$ satisfying $\red_K(0)\le\delta$.
Let $\mathcal J\subset\mathbb N$ be finite, and suppose that for each $j\in\mathcal J$, $K_j$ is a compact subset of
\[
K^\circ\cap\left\{(t,x):G_{b_0}(-t,x)\ge\rho^{-j}\right\}.
\]
Then
\[
\sum_{j\in\mathcal J}\rho^{-j}\capacity(K_j) \le C\red_K(0).
\]
\end{lemma}
\begin{proof}
Choose $\chi\in C_c^\infty(K^\circ)$ with $0 \le \chi \le1$ and $\chi=1$ in a neighborhood of $\bigcup_{j\in\mathcal J} K_j$.
Let $u_n=\mathscr Gf_n$ be as in Lemma~\ref{lem11}. Then
\[
0 \le u_n\le1,\qquad u_n(0)\le\red_K(0)\le\delta,
\]
and
\[
u_n\to1 \quad\text{uniformly on }\;\bigcup_{j\in\mathcal J} K_j.
\]

For $q\ge1$, set $r=q^{-1/d}$.
Since $Pu_n=f_n\ge0$, Duhamel's formula and mass conservation of $\Gamma$ yield
\[
1-\delta\le \int_{\mathbb R^d} \Gamma(0,0,-8r^2,y) \left(1-u_n(-8r^2,y)\right)\,dy.
\]
By the Gaussian upper bound and $0\le1-u_n\le1$,
\[
1-\delta \le C_0\,8^{-d/2}r^{-d} \int_{\abs{y}<\alpha r}\left(1-u_n(-8r^2,y)\right)\,dy+C_0\,8^{-d/2} \int_{\abs{x}\ge \alpha}e^{-a_0\abs{x}^2/8}\,dx,
\]
where we used the change of variables $x=y/r$ in the second integral.
Choose $\alpha\ge 1$ sufficiently large that the second term is at most $(1-\delta)/2$.
Then
\begin{equation}	\label{eq0409sat}
r^{-d}\int_{\abs{y}<\alpha r} \left(1-u_n(-8r^2,y)\right)\,dy\ge c>0.
\end{equation}
Here $\alpha$ and $c$ depend only on $\delta$, $d$, $a_0$, and $C_0$.

Let $\beta$ and $\gamma>0$ be the constants in Lemma~\ref{lem10} corresponding to $\alpha$.
For $q\ge 1$, set
\[
v_q=\mathscr G(\mathbbm1_{\cE_{\beta q}}f_n), \qquad w_q=1-u_n+v_q.
\]
Since $0\le v_q\le u_n\le 1$, we have $0\le w_q\le 1$ and $w_q\ge1-u_n$.
Moreover, $w_q$ is $P$-caloric in $\cE_{\beta q}$.
Thus, by Lemma~\ref{lem10} and \eqref{eq0409sat},
\[
\begin{aligned}
\inf_{\cE_q}w_q &\ge \gamma r^{-d} \int_{\abs{y}<\alpha r}w_q(-8r^2,y)\,dy\\
&\ge \gamma r^{-d} \int_{\abs{y}<\alpha r}(1-u_n(-8r^2,y))\,dy \ge\kappa,
\end{aligned}
\]
where $\kappa=\gamma c>0$ is independent of $q$ and $n$.
Since $w_q$ is continuous in $\cE_{\beta q}$, the same lower bound holds on $\{G_{b_0}(-t,x)\ge q\}$.

Choose $n$ sufficiently large that
\[
u_n\ge1-\kappa/2 \qquad\text{on }\;\bigcup_{j\in\mathcal J}K_j.
\]
Since $K_j\subset \{(t,x): G_{b_0}(-t,x)\ge\rho^{-j}\}$, we obtain
\[
v_{\rho^{-j}}=w_{\rho^{-j}}+u_n-1 \ge\kappa/2 \qquad\text{on }\;K_j,\quad j\in\mathcal J.
\]

Fix $j\in\mathcal J$.
Let $\nu\ge0$ be supported on $K_j$ and satisfy $G_{a_0}^\dagger\nu\le1$ in $\mathbb R^{d+1}$, where $G_{a_0}^\dagger$ is defined in \eqref{eq0205sat}.
Since $v_{\rho^{-j}}\ge\kappa/2$ on $K_j$, Tonelli's theorem gives
\[
\begin{aligned}
\frac{\kappa}{2}\nu(K_j) &\le\int_{K_j}v_{\rho^{-j}}(t,x)\,d\nu(t,x)\\
&=\int_{\cE_{\beta \rho^{-j}}}f_n(s,y) \left(\int_{K_j}\Gamma(t,x,s,y)\,d\nu(t,x)\right) \,dy\,ds.
\end{aligned}
\]
By the Gaussian upper bound,
\[
\begin{aligned}
\frac{\kappa}{2}\nu(K_j) &\le C_0\int_{\cE_{\beta \rho^{-j}}}f_n(s,y) \left(\int_{K_j} G_{a_0}(t-s,x-y)\,d\nu(t,x)\right)\,dy\,ds\\
&=C_0\int_{\cE_{\beta \rho^{-j}}}f_n(s,y) G_{a_0}^\dagger\nu(s,y)\,dy\,ds \le C_0\int_{\cE_{\beta \rho^{-j}}}f_n.
\end{aligned}
\]

Taking the supremum over all such measures $\nu$, we obtain
\[
\mathcal C_{a_0}^\dagger(K_j) \le \frac{2C_0}{\kappa} \int_{\cE_{\beta\rho^{-j}}} f_n.
\]
By Lemma~\ref{lem01}, it follows that
\[
\capacity(K_j) \le C\int_{\cE_{\beta\rho^{-j}} }f_n,
\]
where $C$ is independent of $j$, $n$, and $\mathcal J$.

Summing over $j\in\mathcal J$ and applying Tonelli's theorem, we obtain
\[
\begin{aligned}
\sum_{j\in\mathcal J}\rho^{-j}\capacity(K_j)&\le C\sum_{j\in\mathcal J} \rho^{-j}\int_{\cE_{\beta\rho^{-j}}} f_n\\
&=C\int_{\mathbb R^{d+1}}f_n(t,x) \sum_{\substack{j\in\mathcal J\\ \beta\rho^{-j}<G_{b_0}(-t,x)}} \rho^{-j}\,dx\,dt.
\end{aligned}
\]
Since $0<\rho<1$, the admissible indices $j$ are bounded above.
Summing the geometric series gives
\[
\sum_{\substack{j\in\mathcal J\\ \beta\rho^{-j}<G_{b_0}(-t,x)}} \rho^{-j} \le\frac{G_{b_0}(-t,x)}{\beta(1-\rho)}.
\]
Consequently, by the Gaussian lower bound,
\[
\begin{aligned}
\sum_{j\in\mathcal J}\rho^{-j}\capacity(K_j) &\le C\int_{\mathbb R^{d+1}} G_{b_0}(-t,x)f_n(t,x)\,dx\,dt\\
&\le C\int_{\mathbb R^{d+1}} \Gamma(0,0,t,x)f_n(t,x)\,dx\,dt =Cu_n(0)\le C\red_K(0).
\end{aligned}
\]
\end{proof}

For $j\in\mathbb Z$, define
\begin{equation}	\label{eq1538mon}
F_j:=\overline A_j,\qquad
A_j:=\Omega^\complement\cap\overline Q_{\sqrt\rho} \cap\left\{(t,x):\rho^{-j}\le G_{b_0}(-t,x)\le\rho^{-(j+1)}\right\},
\end{equation}
where
\[
\overline Q_{\sqrt\rho} =\{(t,x):\abs{t}\le\rho,\ \abs{x}\le\sqrt\rho\}.
\]
Since $G_{b_0}$ is continuous away from the origin, $A_j$ is relatively closed in $\mathbb R^{d+1}\setminus\{0\}$.
Thus $F_j$ is compact and
\[
F_j\setminus\{0\}=A_j.
\]

\begin{proposition}	\label{prop03}
If
\[
\sum_{j=1}^\infty \rho^{-j}\capacity(F_j)=\infty,
\]
then $0$ is a regular boundary point of $\Omega$.
\end{proposition}

\begin{proof}
Suppose $0$ is irregular.
By Proposition~\ref{prop02}, there exist $\epsilon>0$ and $0<\delta<1$ such that
\[
\red_{E_\epsilon}(0)<\delta,\qquad E_\epsilon=\Omega^\complement\cap Q_\epsilon^-(0).
\]
By shrinking $\epsilon$ if necessary, we may assume $\epsilon<1$.
Since
\[
\left\{G_{b_0}(-t,x)\ge\rho^{-j}\right\}
\subset
\left\{0<-t\le\rho^{2j/d},\
\abs{x}\le\sqrt{\frac{d}{2eb_0}}\,\rho^{j/d}\right\},
\]
there exists an integer $j_0\ge1$ such that
\[
F_j\setminus\{0\}\subset E_\epsilon
\qquad(j\ge j_0).
\]
Let $N\ge j_0$ be arbitrary.
By Lemma~\ref{lem01}, the thermal capacity of $F_j$ can be approximated by those of compact subsets avoiding the origin:
\[
\capacity(F_j)=\sup\{\capacity(K):K\subset F_j\setminus\{0\},\;K\text{ compact}\}.
\]
For any $\eta>0$, we may therefore choose compact sets $K_j\subset F_j\setminus\{0\}$ such that
\[
\sum_{j=j_0}^N\rho^{-j} \left(\capacity(F_j)-\capacity(K_j)\right)<\eta.
\]

Set $E=\bigcup_{j=j_0}^N K_j$.
By monotonicity,
\[
\red_E(0)\le\red_{E_\epsilon}(0)<\delta.
\]
Since $\red_E(0)<\delta$, Lemma~\ref{lem02} yields a compact neighborhood $K\Subset(-1,0)\times\mathbb R^d$ of $E$ such that $\red_K(0)\le\delta$.
Since
\[
K_j\subset K^\circ\cap \left\{G_{b_0}(-t,x)\ge\rho^{-j}\right\},
\]
Lemma~\ref{lem12} gives
\[
\sum_{j=j_0}^N\rho^{-j}\capacity(K_j) \le C\red_K(0)\le C\delta.
\]
Consequently,
\[
\sum_{j=j_0}^N\rho^{-j}\capacity(F_j) \le C\delta+\eta.
\]
Since $\eta>0$ is arbitrary,
\[
\sum_{j=j_0}^N\rho^{-j}\capacity(F_j) \le C\delta.
\]
The constant $C$ is independent of $N$.
Letting $N\to\infty$ contradicts the divergence of the series.
\end{proof}

\subsection{From Gaussian level sets to the LTU criterion}
By \eqref{eq0237sat}, for any $a>0$,
\begin{equation}	\label{eq1122sun}
\rho^{a\ell -kd/2} \le G_a(-t,x) \le\rho^{-d/2-a}\rho^{a\ell -kd/2} \qquad\text{for }(t,x)\in\Omega_k^\ell.
\end{equation}
Combining this with \eqref{eq1538mon}, we see that $\Omega_k^\ell\cap F_j\neq\emptyset$ implies
\begin{equation}	\label{eq0432sat}
\frac{kd}{2}-b_0\ell-1 \le j \le\frac{kd}{2}-b_0\ell+\frac d2+b_0.
\end{equation}
Consequently, each $\Omega_k^\ell$ intersects only a uniformly bounded number of $F_j$, and for fixed $j,k$, only a uniformly bounded number of $\ell$ can occur.
Moreover, on each nonempty intersection,
\[
\rho^{b_0\ell-kd/2}\asymp\rho^{-j},
\]
with constants independent of the indices.

\begin{proposition}	\label{prop04}
There exists $C>0$, depending only on $d$, $\rho$, and $b_0$, such that
\[
\sum_{k=1}^\infty\sum_{\ell=1}^\infty \capacity(\Omega_k^\ell)\rho^{b_0\ell-kd/2} \le C\capacity(\overline Q_{\sqrt\rho}) +C\sum_{j=1}^\infty\rho^{-j}\capacity(F_j).
\]
\end{proposition}

\begin{proof}
We first show that
\begin{equation}	\label{eq0424sat}
\sum_{k=1}^\infty \capacity(F_j\cap \{\rho^{k+1}\le -t\le\rho^k\}) \le C\capacity(F_j),
\end{equation}
with $C$ independent of $j$.

Fix $j$ and set
\[
\tilde F_j:=\rho^{-j/d}F_j,
\]
where the parabolic scaling is as defined in \eqref{eq1016wed}.
By Gaussian homogeneity,
\begin{equation}	\label{eq1014tue}
1\le G_{b_0}(-t,x)\le\rho^{-1} \qquad ((t,x)\in\tilde F_j\setminus\{0\}).
\end{equation}
Under this scaling, the $k$-th time slab becomes
\[
\{\rho^{k+1-2j/d}\le-t\le\rho^{k-2j/d}\}.
\]
Set $m=\lfloor 2j/d\rfloor$, $i=k-m$, and $\vartheta=\rho^{m-2j/d}$.
Then the time slabs take the form
\[
\vartheta\rho^{i+1}\le-t\le\vartheta\rho^i, \qquad 1\le\vartheta<\rho^{-1}.
\]
By \eqref{eq1014tue} and the definition of $G_{b_0}$, we have $-t\le1$ on $\tilde F_j\setminus\{0\}$.
Thus only the time slabs with $i\ge0$ contribute, except possibly for $i=-1$, whose contribution is a singleton $\{(-1,0)\}$ of zero thermal capacity.

By \eqref{eq1014tue}, for $(t,x)\in\tilde F_j\setminus\{0\}$, we have
\[
\frac{b_0\abs{x}^2}{-t}=\frac{d}{2}\log\frac{1}{-t}-\eta,\qquad\eta=\log G_{b_0}(-t,x)\in [0,\log(1/\rho)].
\]
For $i\ge 0$, set
\[
K_{j,i}:=\tilde F_j\cap\left\{\vartheta\rho^{i+1}\le-t\le\vartheta\rho^i\right\}.
\]
Since $1\le\vartheta<\rho^{-1}$, we have $\rho^{i+1}\le-t\le\rho^{i-1}$ on $K_{j,i}$.
Hence, there exists an integer $i_0\ge1$, depending only on $d$, $b_0$, and $\rho$, such that
\[
c r_i \le \abs{x} \le C r_i,\qquad r_i:=\sqrt{(i+1)\rho^i},
\]
whenever $i\ge i_0$ and $(t,x)\in K_{j,i}$.
Here and below, the constants $c,C>0$ depend only on $d$, $b_0$, and $\rho$ and may vary from line to line.

For each $i\ge0$, let $\mu_i\ge0$ be supported on $K_{j,i}$ and satisfy $G_{b_0}\mu_i\le1$ in $\mathbb R^{d+1}$.
For $i\ge i_0$, since $r_i^2\le r_i^2-t\le Cr_i^2$ and $\abs{x}\le Cr_i$ on $K_{j,i}$, we have
\[
G_{b_0}(r_i^2-t,-x)\ge cr_i^{-d} \qquad ((t,x)\in K_{j,i}).
\]
Evaluating $G_{b_0}\mu_i$ at $(r_i^2,0)$ then gives
\[
\mu_i(\mathbb R^{d+1})\le Cr_i^d.
\]

Choose a sufficiently large integer $M$, depending only on $d$, $b_0$, and $\rho$, such that 
\[
\rho^{M-1}\le \tfrac12 \qquad\text{and}\qquad \abs{x}\le\tfrac12\abs{y}
\]
whenever $(t,x)\in K_{j,k}$ and $(\tau,y)\in K_{j,i}$ with $i\ge i_0$ and $k\ge i+M$.
This is possible since $r_k/r_i\to0$ uniformly as $k-i\to\infty$.
For such points,
\[
\abs{x-y}\ge cr_i,\qquad \tfrac12\rho^{i+1}\le t-\tau\le\rho^{i-1}.
\]
Consequently,
\[
G_{b_0}(t-\tau,x-y) \le C\rho^{-id/2}e^{-c(i+1)}.
\]
Using $\mu_i(\mathbb R^{d+1})\le Cr_i^d$, we obtain
\[
G_{b_0}\mu_i(z) \le C(i+1)^{d/2}e^{-c(i+1)} \qquad (z\in K_{j,k},\ k\ge i+M).
\]

For any finite set of indices $\mathcal I$, let
\[
v=\sum_{i\in \mathcal I} G_{b_0}\mu_i, \qquad S=\bigcup_{i\in \mathcal I} \supp\mu_i.
\]
Fix $\zeta\in S$ belonging to $K_{j,k}$.
The terms with $i<i_0$ or $k-M<i\le k+1$ contribute at most $i_0+M+1$, since $G_{b_0}\mu_i \le 1$.
For $i_0\le i\le k-M$, each $G_{b_0}\mu_i$ is continuous in a neighborhood of $\zeta$, and the preceding estimate gives
\[
\limsup_{z\to\zeta} \sum_{\substack{i\in \mathcal I\\i_0\le i\le k-M}} G_{b_0}\mu_i(z)
\le C\sum_{i=i_0}^\infty (i+1)^{d/2}e^{-c(i+1)} \le C.
\]
The terms with $i\ge k+2$ vanish in a neighborhood of $\zeta$, since their supports lie strictly later in time and $G_{b_0}(s,x)=0$ for $s\le 0$.
Consequently,
\[
\limsup_{z\to\zeta}v(z)\le C \qquad\text{for every }\zeta\in S.
\]
The same estimates also give $v\le C$ on $S$.
Since $v$ satisfies $\partial_t v-(4b_0)^{-1}\Delta v=0$ outside $S$ and $v(t,x)\to0$ as $\abs{x}\to\infty$ uniformly on bounded time intervals, the maximum principle gives $v\le C$ in $\mathbb R^{d+1}$.
Thus
\[
G_{b_0}\left(\frac{1}{C}\sum_{i\in \mathcal I}\mu_i\right) \le1\qquad\text{in }\mathbb R^{d+1}.
\]
Since $C^{-1}\sum_{i\in\mathcal I}\mu_i$ is supported on $\tilde F_j$, the definition of Gaussian capacity gives
\[
\sum_{i\in\mathcal I}\mu_i(K_{j,i}) =\sum_{i\in\mathcal I}\mu_i(\tilde F_j)
\le C\mathcal C_{b_0}(\tilde F_j).
\]
Taking suprema independently over the admissible measures $\mu_i$, we obtain
\[
\sum_{i\in \mathcal I}\mathcal C_{b_0}(K_{j,i}) \le C\mathcal C_{b_0}(\tilde F_j).
\]
Applying Lemma~\ref{lem01}, letting $\mathcal I$ exhaust the nonnegative integers, and undoing the parabolic scaling using \eqref{eq1717mon} yield \eqref{eq0424sat}.

By \eqref{eq0432sat}, there exists a positive integer $N$, depending only on $d$, $b_0$, and $\rho$, such that for each $j,k$, at most $N$ indices $\ell$ satisfy $F_j\cap\Omega_k^\ell\neq\emptyset$.
Moreover, for such indices,
\[
\rho^{b_0\ell-kd/2}\asymp\rho^{-j}.
\]
Since $\Omega_k^\ell\subset \{\rho^{k+1} \le -t \le \rho^k\}$, monotonicity gives
\[
\sum_{\ell=1}^\infty \capacity(F_j\cap\Omega_k^\ell) \le N \capacity(F_j\cap \{\rho^{k+1} \le -t \le \rho^k\}).
\]
Combining these observations with subadditivity and \eqref{eq0424sat}, we obtain
\[
\begin{aligned}
\sum_{k=1}^\infty \sum_{\ell=1}^\infty \rho^{b_0\ell-kd/2} \capacity(\Omega_k^\ell)
&\le C\sum_{j =-\infty}^\infty \rho^{-j} \sum_{k=1}^\infty \sum_{\ell=1}^\infty \capacity(F_j\cap\Omega_k^\ell)\\
&\le C\sum_{j =-\infty}^\infty \rho^{-j}\capacity(F_j).
\end{aligned}
\]
Finally, since $F_j\subset\overline Q_{\sqrt\rho}$,
\[
\sum_{j=-\infty}^0\rho^{-j}\capacity(F_j) \le\frac{1}{1-\rho} \capacity(\overline Q_{\sqrt{\rho}}) <\infty.
\]
This completes the proof.
\end{proof}

\begin{proof}[Proof of Theorem~\ref{thm01}(i)]
By Proposition~\ref{prop04} and the finiteness of $\capacity(\overline Q_{\sqrt\rho})$, the divergence of \eqref{eq:sufficient} implies
\[
\sum_{j=1}^\infty\rho^{-j}\capacity(F_j)=\infty.
\]
Proposition~\ref{prop03} then yields the regularity of $0$.
\end{proof}

\subsection{Necessity}
\begin{proof}[Proof of Theorem~\ref{thm01}(ii)]
For every nonempty $\Omega_k^\ell$, \eqref{eq1122sun} with $a=a_0$ gives
\[
\max_{(t,x)\in\Omega_k^\ell}G_{a_0}(-t,x) \le \rho^{-d/2-a_0}\rho^{a_0\ell-kd/2}.
\]
Hence, Lemma~\ref{lem14} yields
\[
\red_{\Omega_k^\ell}(0) \le C\capacity(\Omega_k^\ell)\rho^{a_0\ell-kd/2}.
\]

Suppose that the series in \eqref{eq:necessary} converges.
For sufficiently small $\epsilon>0$, the set 
\[
E_\epsilon=\Omega^\complement\cap Q_\epsilon^-(0)
\]
is covered by the sets $\Omega_k^\ell$ with $k\ge k(\epsilon)$, where $k(\epsilon)\to\infty$ as $\epsilon\downarrow0$.
By monotonicity and Lemma~\ref{lem02},
\[
\red_{E_\epsilon}(0) \le \sum_{k\ge k(\epsilon)}
\sum_{\ell=1}^\infty
\red_{\Omega_k^\ell}(0) \le C\sum_{k\ge k(\epsilon)}\sum_{\ell=1}^\infty \capacity(\Omega_k^\ell)\rho^{a_0\ell-kd/2} \to0.
\]
This contradicts Proposition~\ref{prop02}, since $0$ is regular.
Therefore, \eqref{eq:necessary} holds.
\end{proof}
\subsection{Boundary regularity in time-varying domains}
We first establish the behavior of thermal capacity under the transformation $(t,x)\mapsto (t, \Phi(t,x))$.

\begin{lemma}	\label{lem15}
Let $I\subset\mathbb R$ be an open interval and $B\subset\mathbb R^d$ be a ball.
Suppose that $\Phi:I\times B\to\mathbb R^d$ satisfies \eqref{eq1254mon} for some $N\ge 1$ and all $t, s\in I$ and $x, y \in B$.
Set $\mathsf T(t,x)=(t,\Phi(t,x))$.
Then
\[
\capacity(\mathsf T(K))\asymp\capacity(K)
\]
for every compact $K\subset I\times B$, with constants depending only on $d$ and $N$.
\end{lemma}

\begin{proof}
For $z=(t,x)$ and $\zeta=(s,y)$ in $I\times B$ with $t>s$, \eqref{eq1254mon} gives
\[
N^{-1}\abs{x-y}-N\sqrt{t-s} \le \abs{\Phi(t,x)-\Phi(s,y)} \le N\abs{x-y}+N\sqrt{t-s}.
\]
Using $(a+b)^2\le2a^2+2b^2$, we obtain
\begin{equation}	\label{eq1310mon}
e^{-N^2/2}G_{N^2/2}(z-\zeta) \le G_{1/4}(\mathsf Tz-\mathsf T\zeta) \le e^{N^2/4}G_{1/(8N^2)}(z-\zeta).
\end{equation}
The upper bound on $\abs{\Phi(t,x)-\Phi(s,y)}$ also gives
\begin{equation}	\label{eq1311mon}
G_{1/4}(z-\zeta) \le e^{1/4}G_{1/(8N^2)}(\mathsf Tz-\mathsf T\zeta).
\end{equation}

Let $\mu\ge0$ be supported on $K$ with $G_{1/(8N^2)}\mu\le1$, and set $\nu=\mathsf T_\#\mu$.
By \eqref{eq1310mon},
\[
G_{1/4}\nu\le e^{N^2/4} \quad\text{in }\,\mathsf T(I\times B).
\]
Since $\supp\nu\Subset \mathsf T(I\times B)$, the potential $G_{1/4}\nu$ is caloric outside $\mathsf T(I\times B)$ and vanishes as $\abs{x}\to\infty$.
Thus the maximum principle extends this bound to $\mathbb R^{d+1}$.
Hence, by Lemma~\ref{lem01},
\[
\capacity(\mathsf T(K)) \ge c\mathcal C_{1/(8N^2)}(K) \ge c\capacity(K).
\]

Conversely, let $\nu\ge0$ be supported on $\mathsf T(K)$ with $G_{1/(8N^2)}\nu\le1$, and set $\mu=(\mathsf T^{-1})_\#\nu$.
By \eqref{eq1311mon} and the same argument, $G_{1/4}\mu\le e^{1/4}$ on $\mathbb R^{d+1}$.
Hence
\[
\capacity(K) \ge c\mathcal C_{1/(8N^2)}(\mathsf T(K)) \ge c\capacity(\mathsf T(K)).\qedhere
\]
\end{proof}

For an open set $\Omega\subset\mathbb R^{d+1}$, $z\in\partial \Omega$, and $b>0$, define the LTU series
\[
\mathscr W^{(b)}(\Omega,z)=\sum_{k=1}^\infty\sum_{\ell=1}^\infty \capacity(\Omega_k^\ell(z,\rho))\,\rho^{b\ell-kd/2},
\]
where $\Omega_k^\ell(z,\rho)$ is defined in \eqref{eq0237sat} and $0<\rho<1$ is fixed.

\begin{proof}[Proof of Corollary~\ref{cor02}]
Since boundary regularity is a local property by Proposition~\ref{prop02}, we may assume that $\Omega$ is bounded.
Set $\zeta_0=(t_0,y_0)$ and choose $\epsilon>0$ sufficiently small that
\[
\mathcal C=(t_0-\epsilon^2,t_0+\epsilon^2)\times(D\cap B_\epsilon(y_0))
\]
is contained in $U$ and
\[
\mathsf T(Q_\epsilon^-(\zeta_0))\subset Q_\delta^-(z_0).
\]
We denote the LTU sets associated with $\mathcal C$ by $\mathcal C_k^m(\zeta_0,\rho)$.

Since $\mathsf T$ is a homeomorphism onto its open image, we may choose $\eta>0$ sufficiently small that
\[
Q_\eta^-(z_0)\subset \mathsf T(Q_\epsilon^-(\zeta_0)), \qquad \mathsf T(Q_\eta^-(\zeta_0))\subset Q_\delta^-(z_0).
\]

Restricting the LTU sets to $\overline Q_\eta(z_0)$ and $\overline Q_\eta(\zeta_0)$, respectively, changes each LTU series by at most a finite quantity.

We claim that for every $b>0$, there exist $C\ge 1$, independent of $b$, and $C_b>0$ such that
\begin{equation}	\label{eq1344mon}
\begin{aligned}
\mathscr W^{(b)}(\Omega,z_0) &\le C_b\bigl(1+\mathscr W^{(b/C)}(\mathcal C,\zeta_0)\bigr),\\
\mathscr W^{(b)}(\mathcal C,\zeta_0) &\le C_b\bigl(1+\mathscr W^{(b/C)}(\Omega,z_0)\bigr).
\end{aligned}
\end{equation}
We also take $\eta$ sufficiently small that the localized sets and their images and preimages satisfy the distance cutoff in \eqref{eq0237sat}.

By Lemma~\ref{lem15}, thermal capacities are comparable under $\mathsf T$ in these neighborhoods.
For $s<t_0$ and $y$ sufficiently close to $y_0$, \eqref{eq1254mon} gives
\[
\begin{aligned}
\frac{\abs{y-y_0}^2}{t_0-s} &\le 2N^2 \frac{\abs{\Phi(s,y)-x_0}^2}{t_0-s}+2N^4,\\
\frac{\abs{\Phi(s,y)-x_0}^2}{t_0-s} &\le 2N^2 \frac{\abs{y-y_0}^2}{t_0-s}+2N^2.
\end{aligned}
\]

Since $\mathsf T$ preserves the time coordinate, the time index $k$ is unchanged. Moreover, the first estimate gives
\[
\frac{\abs{y-y_0}^2}{t_0-s} \le 2N^2\ell\log(1/\rho)+2N^4
\]
whenever $\mathsf T(s,y)\in\Omega_k^\ell(z_0,\rho)$.
Thus,
\[
\mathsf T^{-1}\left(\Omega_k^\ell(z_0,\rho)\cap\overline Q_\eta(z_0)\right)\subset\bigcup_{1\le m\le C(\ell+1)}\mathcal C_k^m(\zeta_0,\rho).
\]
Here the local identity of the two domains ensures that the preimage belongs to $\mathcal C^\complement$.
By Lemma~\ref{lem15} and subadditivity of capacity,
\[
\capacity\left(\Omega_k^\ell(z_0,\rho)\cap\overline Q_\eta(z_0)\right)\le C\sum_{m\le C(\ell+1)}\capacity\left(\mathcal C_k^m(\zeta_0,\rho)\right).
\]
Summing over $\ell$ and using
\[
\sum_{\ell\ge m/C-1}\rho^{b\ell}\le C_b\rho^{bm/C},
\]
we obtain the first inequality in \eqref{eq1344mon}.

Similarly, the second spatial estimate gives
\[
\mathsf T\left(\mathcal C_k^m(\zeta_0,\rho)\cap\overline Q_\eta(\zeta_0)\right)
\subset\bigcup_{1\le\ell\le C(m+1)}\Omega_k^\ell(z_0,\rho).
\]
The second inequality in \eqref{eq1344mon} follows by the same argument.

For every $b>0$, the operator $\partial_t-(4b)^{-1}\Delta$ has Gaussian exponents $a_0=b_0=b$.
The pointwise version of Tikhonov's classical theorem \cite{Tyc1938} asserts that regularity of a lateral boundary point of a cylindrical domain for the heat operator is equivalent to regularity of the corresponding spatial boundary point for the Laplacian.
Thus, applying Theorem~\ref{thm01} to this operator in $\mathcal C$, we obtain
\[
\mathscr W^{(b)}(\mathcal C,\zeta_0)=\infty \iff y_0\text{ is regular for }\Delta\text{ in }D.
\]

If $y_0$ is regular for $\Delta$, then $\mathscr W^{(Cb_0)}(\mathcal C,\zeta_0)=\infty$.
By \eqref{eq1344mon}, $\mathscr W^{(b_0)}(\Omega,z_0)=\infty$, so $z_0$ is regular for $P$ by Theorem~\ref{thm01}(i).

Conversely, if $y_0$ is irregular for $\Delta$, then $\mathscr W^{(a_0/C)}(\mathcal C,\zeta_0)<\infty$.
By \eqref{eq1344mon}, $\mathscr W^{(a_0)}(\Omega,z_0)<\infty$, so $z_0$ is irregular for $P$ by Theorem~\ref{thm01}(ii).
The assertion for the heat operator follows by taking $P=\partial_t-\Delta$.
\end{proof}

\section{A quantitative Wiener-type criterion}	\label{sec05}
Adapting the capacitary packing and stopping-time arguments recently developed in \cite{Kim2026} for the heat equation, we establish, under the hypotheses of Theorem~\ref{thm01}, a quantitative boundary estimate in terms of finite partial sums of the Wiener-type series at the critical exponent $b_0$.

Fix $0<\rho<1$ and, for $z_0\in\partial\Omega$, set
\[
\mathscr W_n(z_0)=\sum_{k=1}^n\sum_{\ell=1}^\infty \capacity(\Omega_k^\ell(z_0,\rho))\rho^{b_0\ell-kd/2},
\qquad \mathscr W_n(z_0)=0\quad(n \le 0).
\]
For $f\in C(\partial\Omega)$ and $r>0$, define
\[
\omega_f(z_0,r)=\sup_{\zeta\in\partial\Omega\cap\overline Q_r(z_0)} \abs{f(\zeta)-f(z_0)},
\]
and, for $z\ne z_0$, set
\[
N(z)=\max\left\{0,\left\lfloor\frac{2\log(1/\abs{z-z_0})}{\log(1/\rho)}\right\rfloor\right\}.
\]

\begin{theorem}	\label{thm02}
Let $P$ and $\Omega$ be as in Theorem~\ref{thm01}, $z_0\in\partial\Omega$, and $f\in C(\partial\Omega)$.
There exist constants $C, c>0$, depending only on $d$, $\rho$, and the constants in Condition~\ref{gaussian}, such that, for $z\in\Omega$ and integers $m\ge 1$,
\begin{equation}	\label{eq0528mon}
\abs{\overline H_f(z)-f(z_0)} \le C\left\{\omega_f(z_0,\rho^{m/2})+\norm{f}_\infty\exp\left(-c(\mathscr W_{N(z)}(z_0)-\mathscr W_{m-1}(z_0))_+\right)\right\}.
\end{equation}
In particular,
\begin{equation}\label{eq1121tue}
\abs{\overline H_f(z)-f(z_0)}\le C\left\{\omega_f\left(z_0,e^{-c\mathscr W_{N(z)}(z_0)}\right)+\norm{f}_\infty e^{-c\mathscr W_{N(z)}(z_0)}\right\}.
\end{equation}
The same estimates hold for $\underline H_f$.
\end{theorem}

The estimate \eqref{eq1121tue} also yields Theorem~\ref{thm01}(i).

\subsection{An exponential capacitary estimate}
By translation, we may assume $z_0=0$.
Let $F_j$ be as in \eqref{eq1538mon} and set
\[
r_j=\rho^{j/d},\qquad s_j=-8r_j^2.
\]
Let $\beta$ and $\gamma$ be the constants in Lemma~\ref{lem10} corresponding to $\alpha=1$.
Set
\[
\sigma=1+\max\{16,\beta^{-2/d}\}.
\]
In the following proof, we use the Gaussian bounds for $0<t-\tau<2\sigma$, adjusting $c_0$ and $C_0$ if necessary.
For $\varpi>2$ to be fixed below, put
\[
\mathcal S_j=(-\sigma r_j^2,\sigma r_j^2)\times B_{\varpi r_j}(0).
\]

\begin{lemma}		\label{lem16}
There exist $\varpi>2$, an integer $\ell_0 \ge1$, and $C,c>0$, depending only on $d$, $\rho$, and the constants in Condition~\ref{gaussian}, such that for every compact set $K\subset\mathbb R^{d+1}$ and integers $1\le m\le n$,
\begin{equation}	\label{eq0631mon}
\sup_{\mathcal S_{n+\ell_0}}(1-\bal_K)
\le C\exp\left(-c\sum_{j=m}^n \rho^{-j}\capacity(K\cap F_j)\right).
\end{equation}
\end{lemma}

\begin{proof}
We first work with compact sets $E_j\subset (K\cap F_j)\setminus\{0\}$, $m\le j\le n$, and functions of the form $v=1-\mathscr Gg$, where $g\ge 0$ is bounded and compactly supported and $0\le v\le1$ in $\mathbb R^{d+1}$.
Let $\epsilon\in(0,1)$ be arbitrary and suppose that $v\le\epsilon$ on $\bigcup_{j=m}^n E_j$.
Set
\[
\alpha_j=\fint_{B_{r_j}}v(s_j,y)\,dy,\qquad M_j=\sup_{\mathcal S_j}v\qquad(j\ge1),
\]
and
\[
\kappa_j=\rho^{-j}\capacity(E_j) \qquad (m\le j\le n).
\]
The Gaussian level-set geometry and \eqref{eq1717mon} imply that $\kappa_j\le A_0$ for a constant $A_0$ independent of $j$ and $E_j$.

\subsubsection*{Step 1. Capacitary packing.}
Set $w_j=\mathscr G(\mathbbm1_{\cE_{\beta\rho^{-j}}}g)$ and $h_j=v+w_j$.
Since
\[
0 \le w_j\le\mathscr Gg=1-v,
\]
the function $h_j$ is bounded and nonnegative.
Moreover, $Ph_j=0$ in $\cE_{\beta\rho^{-j}}$.
Lemma~\ref{lem10} and continuity give $h_j\ge\gamma\abs{B_1}\alpha_j$ on $E_j$.
Hence
\[
w_j\ge\gamma\abs{B_1}\alpha_j-\epsilon \quad\text{on }\;E_j.
\]
As in the proof of Lemma~\ref{lem12}, integration against measures $\nu$ supported on $E_j$ satisfying $G_{a_0}^\dagger\nu\le1$, followed by Lemma~\ref{lem01}, yields
\[
\alpha_j\capacity(E_j)\le C\int_{\cE_{\beta\rho^{-j}}}g+C\epsilon\capacity(E_j).
\]
Indeed, Tonelli's theorem and the Gaussian upper bound give
\[
\gamma\abs{B_1}\alpha_j\nu(E_j) \le C_0\int_{\cE_{\beta\rho^{-j}}}g+\epsilon\nu(E_j),
\]
and we may take a sequence of admissible measures $\nu$ such that $\nu(E_j)\to\mathcal C_{a_0}^\dagger(E_j)$.
Since $\cE_{\beta\rho^{-j}}=\{G_{b_0}(-t,x)>\beta\rho^{-j}\}$, the geometric series estimate
\[
\sum_{j=m}^n \rho^{-j}\mathbbm1_{\cE_{\beta\rho^{-j}}}(t,x) \le\frac{G_{b_0}(-t,x)}{\beta(1-\rho)}\mathbbm1_{\cE_{\beta\rho^{-m}}}(t,x)
\]
and the Gaussian lower bound yield
\begin{equation}	\label{eq0703mon}
\sum_{j=m}^n\alpha_j\kappa_j\le Cw_m(0)+C\epsilon\sum_{j=m}^n\kappa_j.
\end{equation}

We claim that $w_m(0)\le2M_m$ if $\varpi$ is sufficiently large.
In fact,
\[
\cE_{\beta\rho^{-m}}\subset(-\beta^{-2/d}r_m^2,0)\times B_{\varpi_0r_m}, \qquad \varpi_0=\beta^{-1/d}\sqrt{\frac{d}{2eb_0}},
\]
and
\[
w_m(0)\ge c_0\beta\rho^{-m}\int_{\cE_{\beta\rho^{-m}}}g.
\]
Consider $Q=(-\sigma r_m^2,0)\times B_{\varpi r_m}$.
The function $w_m$ vanishes on the bottom of $Q$.
For $\varpi>2\varpi_0$, the bound
$G_{a_0}(t,x)\le C\abs{x}^{-d}$ for $x\ne0$ yields
\[
w_m\le C\varpi^{-d}\rho^{-m} \int_{\cE_{\beta\rho^{-m}}}g \le C\varpi^{-d}w_m(0)
\]
on the lateral boundary of $Q$.
Choose $\varpi$ large enough that the last expression is at most $w_m(0)/2$.
Since
\[
Ph_m=-\mathbbm1_{\cE_{\beta\rho^{-m}}^\complement}g \le 0,
\]
comparison and continuity at $0$ give
\[
w_m(0)\le h_m(0)\le M_m+\tfrac12 w_m(0).
\]
Combining this with \eqref{eq0703mon}, we obtain
\begin{equation}	\label{eq0733mon}
\sum_{j=m}^n\alpha_j\kappa_j\le C M_m+C\epsilon\sum_{j=m}^n\kappa_j.
\end{equation}
The constants depend only on $d$, $\rho$, and the constants in Condition~\ref{gaussian}.
The same estimate holds with $m, n$ replaced by any integers $m',n'$ satisfying $m\le m'\le n'\le n$.

\subsubsection*{Step 2. Contraction and stopping-time argument.}
We show that, after increasing $\varpi$ and choosing $\ell_0$ sufficiently large, there exists $\eta\in(0,1)$ such that
\begin{equation}	\label{eq0737mon}
M_{j+\ell_0}\le\eta M_j+(1-\eta)\alpha_j,\qquad j\ge 1.
\end{equation}
Choose $\ell_0$ so that $\sigma\rho^{2\ell_0/d}<4$ and $\varpi\rho^{\ell_0/d}<1$.
Choose a continuous cutoff $\phi_j$ supported in $B_{2r_j}$, with $0\le\phi_j\le1$ and $\phi_j=1$ on $B_{r_j}$.
Put
\[
\psi_j=(M_j-v(s_j,\cdot))\phi_j \quad\text{on }B_{2r_j}
\]
and extend it by zero.
By the Gaussian lower bound, for $(t,x)\in\mathcal S_{j+\ell_0}$,
\[
\mathscr S_{s_j}\psi_j(t,x)\ge c_1 \rho^{-j}\int\psi_j,
\]
where $c_1=c_0 12^{-d/2}\exp(-9b_0/4)$, since $4r_j^2<t-s_j<12r_j^2$ and $\abs{x-y}<3r_j$ for $y\in\supp\psi_j$.
On the lateral boundary of $(s_j,4r_j^2)\times B_{\varpi r_j}$, the Gaussian upper bound gives
\[
\mathscr S_{s_j}\psi_j\le C(\varpi-2)^{-d}\rho^{-j}\int\psi_j.
\]
Increase $\varpi$ so that $C(\varpi-2)^{-d}\le c_1/2$, and choose $\ell_0$ as above for this value of $\varpi$.
Since $M_j-v$ is nonnegative and $P$-superparabolic in $\mathcal S_j$, comparison in $(s_j,4r_j^2)\times B_{\varpi r_j}$ gives
\[
M_j-v\ge\mathscr S_{s_j}\psi_j-\frac{c_1}{2}\rho^{-j}\int\psi_j.
\]
In particular, on $\mathcal S_{j+\ell_0}$,
\[
M_j-v\ge \frac{c_1}{2}\rho^{-j}\int\psi_j \ge\frac{c_1}{2}\abs{B_1} (M_j-\alpha_j).
\]
Reducing the constant if necessary proves \eqref{eq0737mon}.

Fix $A_\ast\ge 2C$, where $C$ is the first constant on the right-hand side of \eqref{eq0733mon}.
If
\[
\sum_{j=m}^n\kappa_j\ge A_\ast,
\]
then \eqref{eq0733mon} gives an index $j\in\{m,\ldots,n\}$ such that $\alpha_j \le M_m/2+C\epsilon$.
By \eqref{eq0737mon} and the fact that $(M_j)$ is nonincreasing,
\begin{equation}	\label{eq0747mon}
M_{n+\ell_0}\le\theta M_m+C\epsilon,\qquad \theta=(1+\eta)/2<1.
\end{equation}
Set $m_1=m$ and define recursively
\[
n_i=\min\left\{k\in\{m_i,\ldots,n\}:\sum_{j=m_i}^k\kappa_j\ge A_\ast\right\},\qquad m_{i+1}=n_i+\ell_0,
\]
as long as the defining set is nonempty.
If $N$ is the number of blocks so obtained, then by the minimality of $n_i$, each block, together with the subsequent skipped indices, contributes at most $A_\ast+\ell_0 A_0$.
The remaining sum is less than $A_\ast$.
Consequently,
\[
\sum_{j=m}^n\kappa_j\le N(A_\ast+\ell_0 A_0)+A_\ast.
\]
Iterating \eqref{eq0747mon} and using $m_{i+1}=n_i+\ell_0$, we obtain
\begin{equation}	\label{eq0757mon}
M_{n+\ell_0}\le\theta^N M_m+\frac{C\epsilon}{1-\theta} \le C M_m\exp\left(-c\sum_{j=m}^n\kappa_j\right)+C\epsilon.
\end{equation}
The estimate remains valid when $N=0$.

\subsubsection*{Step 3. Passage to balayage.}
Set $E=\bigcup_{j=m}^n E_j$.
If $E$ is empty, the desired estimate is immediate.
Choose a compact neighborhood $K'\Subset(-\infty,0)\times\mathbb R^d$ of $E$ and a cutoff $\chi$ as in Lemma~\ref{lem11}, with $\chi=1$ near $E$.
Let $u_k=\mathscr Gf_k$ be the approximations in that lemma, and apply \eqref{eq0757mon} with $g=f_k$, $v=1-u_k$, and
\[
\epsilon_k=\sup_E(1-u_k).
\]
Then $0\le v\le 1$ and $\epsilon_k \to 0$ as $k \to \infty$.
Moreover, Lemma~\ref{lem11} gives $u_k\le\red_{K'}$.
Thus, letting $k\to\infty$ in \eqref{eq0757mon}, we obtain
\[
1-\red_{K'}\le C\exp\left(-c\sum_{j=m}^n \rho^{-j}\capacity(E_j)\right)\qquad\text{in }\,\mathcal S_{n+\ell_0}.
\]
Letting $K'\downarrow E$ and using Lemma~\ref{lem02}, we obtain the same bound for $1-\red_E$ and hence, by monotonicity, for $1-\red_K$.
Since $\mathcal S_{n+\ell_0}$ is open, the bound also holds for $1-\bal_K$.
Finally, by Lemma~\ref{lem01}, the capacities of $E_j$ can be taken arbitrarily close to those of $K\cap F_j$, since removing $Q_\delta(0)$ reduces the capacity by at most $C\delta^d$.
This proves \eqref{eq0631mon}.
\end{proof}

\subsection{The boundary modulus of continuity}
For simplicity, write $\mathscr W_n=\mathscr W_n(0)$ and $\Omega_k^\ell=\Omega_k^\ell(0,\rho)$.
We first note that
\begin{equation}	\label{eq0804mon}
\capacity(\Omega_k^\ell)\rho^{b_0\ell-kd/2} \le C(1+\ell)^{d/2}\rho^{b_0\ell}, \qquad \mathscr W_n-\mathscr W_{n-1}\le C.
\end{equation}
Indeed, $\Omega_k^\ell\subset Q_{C\sqrt{1+\ell}\,\rho^{k/2}}(0)$, so the first inequality follows from scaling and monotonicity of capacity, and the second by summing over $\ell$.

For $m\ge 1$, set $K_m=\Omega^\complement\cap\overline Q_{\rho^{m/2}}(0)$.
We claim that
\begin{equation}	\label{eq0815mon}
1-\bal_{K_m}(z)\le C\exp\left(-c(\mathscr W_{N(z)}-\mathscr W_{m-1})_+\right),\qquad z\ne0.
\end{equation}
It suffices to consider $N=N(z)\ge m$.
For $k\ge m$, if $\Omega_k^\ell\cap\{\abs{x}\ge\rho^{m/2}\}\ne\emptyset$, then
\[
\ell\ge\frac{\rho^{m-k}}{\log(1/\rho)}.
\]
By monotonicity of capacity and \eqref{eq0804mon},
\[
\sum_{k=m}^N\sum_{\ell=1}^\infty \capacity(\Omega_k^\ell\cap\{\abs{x}\ge\rho^{m/2}\}) \rho^{b_0\ell-kd/2} \le C\sum_{h=0}^\infty \sum_{\ell\log(1/\rho)\ge\rho^{-h}}(1+\ell)^{d/2}\rho^{b_0\ell}<\infty.
\]
Since $-t\le\rho^m$ on $\Omega_k^\ell$ for $k\ge m$, the subadditivity of capacity gives
\begin{equation}	\label{eq0826mon}
\mathscr W_N-\mathscr W_{m-1} \le C+\sum_{k=m}^N\sum_{\ell=1}^\infty \capacity(K_m\cap\Omega_k^\ell)\rho^{b_0\ell-kd/2}.
\end{equation}

The proof of \eqref{eq0424sat} applies without change to any compact subset of $F_j$.
Combining this observation with \eqref{eq0432sat} and the argument in Proposition~\ref{prop04}, we obtain the finite-scale inequality
\begin{equation}	\label{eq1339mon}
\sum_{k=m}^N\sum_{\ell=1}^\infty\capacity(K_m\cap\Omega_k^\ell)\rho^{b_0\ell-kd/2}
\le C\sum_{j\le j_N}\rho^{-j}\capacity(K_m\cap F_j),
\end{equation}
where $j_N=\lfloor d(N+1)/2\rfloor$.
Indeed, only $j\le j_N$ can occur in intersections with $\Omega_k^\ell$ for $k\le N$, and the estimates of Proposition~\ref{prop04} apply to $K_m\cap F_j$ in place of $F_j$.

Set $n=\lfloor dN/2\rfloor-\ell_0$.
In view of \eqref{eq0804mon} and $1-\bal_{K_m}\le1$, we may assume that $n\ge1$.
The contribution of $j\le0$ in \eqref{eq1339mon} is bounded by $\capacity(K_m)\sum_{j\le0}\rho^{-j}\le C$.
The contribution of $n<j\le j_N$ is also bounded, since $j_N-n\le C$ and $\rho^{-j}\capacity(F_j)\le C$.
It follows from \eqref{eq0826mon} and \eqref{eq1339mon} that
\[
\mathscr W_N-\mathscr W_{m-1} \le C+C\sum_{j=1}^n\rho^{-j}\capacity(K_m\cap F_j).
\]
Since $\abs{z}\le\rho^{N/2}\le\rho^{(n+\ell_0)/d}$, we have $z\in\mathcal S_{n+\ell_0}$.
Thus \eqref{eq0815mon} follows from Lemma~\ref{lem16}.

\begin{proof}[Proof of Theorem~\ref{thm02}]
Let $M$ be as in the construction of $w_r$ in the proof of Proposition~\ref{prop01}.
Choose an integer $h\ge1$ such that $M\rho^{h/2}\le 1$.
For $m\ge1$, set
\[
r=\rho^{(m+h)/2}.
\]
By \eqref{eq1503sun}, $\bal_{K_{m+h}}\le w_r$ and $w_r\le1/4$ on $\partial Q_{Mr}$.
Set
\[
v_r=1-\bal_{K_{m+h}}.
\]
Then $v_r$ is nonnegative and $P$-caloric in $\Omega$, and $v_r\ge1/2$ in a neighborhood of $\partial Q_{Mr}$ by continuity of $w_r$.
Set $A=2\norm{f}_\infty$ and define
\[
w(z)=\left\{
\begin{array}{cl}
\min\left\{A,\omega_f(0,\rho^{m/2})+2Av_r(z)\right\}, &z\in\Omega\cap Q_{Mr},\\
A,&z\in\Omega\setminus Q_{Mr}.
\end{array}
\right.
\]
As in the proof of Proposition~\ref{prop01}, the two expressions agree near $\Omega\cap\partial Q_{Mr}$, so $w$ is $P$-superparabolic in $\Omega$.
Since $Mr\le\rho^{m/2}$, we also have
\[
\liminf_{\Omega\ni z\to\zeta}w(z)\ge\abs{f(\zeta)-f(0)},\qquad\zeta\in\partial\Omega.
\]
Therefore $w$ belongs to the upper Perron classes for both $f-f(0)$ and $f(0)-f$, and bounds $\abs{\overline H_f-f(0)}$ and $\abs{\underline H_f-f(0)}$ from above.

For $z\in\Omega\cap Q_{Mr}$, \eqref{eq0815mon} with $m+h$ in place of $m$ gives the desired bound \eqref{eq0528mon}, since by \eqref{eq0804mon}, $\mathscr W_{m+h-1}-\mathscr W_{m-1}\le Ch$.
For $z\in\Omega\setminus Q_{Mr}$, we have $\abs{z}\ge Mr>\rho^{(m+h)/2}$, so $N(z)<m+h$ and $(\mathscr W_{N(z)}-\mathscr W_{m-1})_+\le Ch$.
Thus \eqref{eq0528mon} follows in this case from the bound $A$.

Finally, by \eqref{eq0804mon}, there exists $C_1>0$ such that
\[
\exp\left(-c(\mathscr W_{N(z)}-\mathscr W_{m-1})_+\right) \le\exp(C_1(m-1)-c\mathscr W_{N(z)}).
\]
Choose $m=\max\{1,\lceil c\mathscr W_{N(z)}/(2C_1)\rceil\}$.
Then
\[
\rho^{m/2}\le e^{-c\log(1/\rho)\mathscr W_{N(z)}/(4C_1)},\qquad C_1(m-1)-c\mathscr W_{N(z)}\le-\frac{c}{2}\mathscr W_{N(z)}.
\]
Substitution into \eqref{eq0528mon} proves \eqref{eq1121tue} with a possibly smaller $c$.
\end{proof}

\appendix
\section{Proofs of preliminary lemmas}\label{sec_a}

\begin{proof}[Proof of Lemma~\ref{lem03}]
(a) It suffices to prove the second assertion.
Suppose that $v$ is negative somewhere in $D$.
For sufficiently small $\epsilon>0$, the function
\[
v_\epsilon(t,x)=v(t,x)+\frac{\epsilon}{T-t}
\]
is negative somewhere in $D$.
Since $v$ is bounded below and satisfies the assumed boundary condition, $v_\epsilon$ attains a negative minimum $m$ at some $z^\circ\in D$.
Set $g(t)=m-\epsilon/(T-t)$.
Then $v\ge g$ in $D$, with equality at $z^\circ$, and $Pg=-\epsilon/(T-t)^2<0$.
Choose a cylinder $Q\Subset D$ containing $z^\circ$ and let $h$ solve $Ph=0$ in $Q$ with boundary data $g$.
By the strong minimum principle, $h>g$ in $Q$.
On the other hand, since $v\ge h$ on $\partial_pQ$, we have $v\ge h$ in $Q$.
Thus $v(z^\circ)\ge h(z^\circ)>g(z^\circ)=v(z^\circ)$, a contradiction.

(b) The assertions for positive scalar multiples and minima follow directly from the definition.
To prove the assertion for sums, let $h\in C(\overline Q)$ be $P$-caloric in a cylinder $Q\Subset D$ with $h\le u+v$ on $\partial_p Q$.
Then
\[
h=\mathscr H_Q h \le\mathscr H_Q u+\mathscr H_Q v \le u+v \qquad\text{in }\,Q.
\]
Since $u+v$ is lower semicontinuous, it is $P$-superparabolic.

(c) The increasing limit $v$ is lower semicontinuous.
Let $h\in C(\overline Q)$ be $P$-caloric in a cylinder $Q\Subset D$ with $h\le v$ on $\partial_p Q$.
For each $\epsilon>0$, lower semicontinuity and compactness of $\partial_p Q$ give an $N$ such that $v_N \ge h-\epsilon$ on $\partial_p Q$.
By comparison, $v\ge v_N \ge h-\epsilon$ in $Q$.
Letting $\epsilon\downarrow 0$ proves the assertion.

(d) Let $Q=(s,t)\times B\Subset D$ and choose $\delta>0$ such that $Q_+=(s,t+\delta)\times B\Subset D$.
Choose continuous functions $\phi_n$ on $\overline Q_+$ increasing pointwise to $v$.
Let $h_n=\mathscr H_{Q_+}\phi_n$.
By comparison, $h_n\le h_{n+1}\le v$ in $Q_+$.
The Harnack convergence theorem implies that $h=\lim_n h_n$ is $P$-caloric in $Q_+$.
Since $h_n=\mathscr H_Q\phi_n$ in $Q$, monotone convergence gives $h=\mathscr H_Qv$ in $Q$.
Moreover, continuity of $h_n$ up to the parabolic boundary and $h \ge h_n$ yield
\begin{equation}	\label{eq1621sun}
\liminf_{Q\ni z\to\zeta}h(z)\ge v(\zeta) \qquad(\zeta\in\partial_pQ).
\end{equation}
Together with $h\le v$ in $Q_+$, this shows that $v^Q \le v$ and that $v^Q$ is lower semicontinuous.

To verify the comparison property, let $w\in C(\overline{Q'})$ be $P$-caloric in a cylinder $Q'\Subset D$ with $w\le v^Q$ on $\partial_pQ'$.
Since $v^Q\le v$, we have $w\le v$ in $Q'$.
If $Q\cap Q'$ is nonempty, then $h-w$ has nonnegative lower boundary limits on its initial and lateral boundaries by \eqref{eq1621sun} and the boundary condition on $w$.
Part (a) gives $w\le h$ in $Q\cap Q'$ and hence, by continuity, on $Q'\cap(\{t\}\times B)$.
Elsewhere in $Q'$, we have $w\le v=v^Q$.
Thus $w\le v^Q$ in $Q'$, as desired.
\end{proof}

\begin{proof}[Proof of Lemma~\ref{lem04}]
Set $u=\inf_{v\in\mathcal F}v$.
Fix a cylinder $Q\Subset D$ and a countable dense subset $\{z_j\}_{j\ge1}$ of $Q$.
Since $\mathcal F$ is closed under finite minima, there exists a decreasing sequence $v_n\in\mathcal F$ such that
\[
u(z_j)\le v_n(z_j)\le u(z_j)+1/n,\qquad 1\le j\le n.
\]
Let $h_n=v_n^Q$.
Then $h_n\in\mathcal F$ is $P$-caloric in $Q$ and
\[
u\le h_{n+1}\le h_n\le v_n \qquad\text{in }\,Q.
\]
By the Harnack convergence theorem, $h_n$ converges locally uniformly to a $P$-caloric function $h$ in $Q$, with $h(z_j)=u(z_j)$ for every $j$.

For any $v\in\mathcal F$, let $w_n=(\min\{v_n,v\})^Q$.
Then $w_n$ decreases locally uniformly to a $P$-caloric function $w$ satisfying $u \le w \le h$ and $w \le v$ in $Q$.
Since $w(z_j)=h(z_j)$ for every $j$, we have $w=h$ in $Q$.
Thus $h \le v$ for every $v\in\mathcal F$, and hence $h=u$ in $Q$.
Since $Q$ was arbitrary, $u$ is $P$-caloric in $D$.
\end{proof}

\begin{proof}[Proof of Lemma~\ref{lem01}]
The equivalence $\mathcal C_a(E)\asymp\capacity(E)$ follows from the comparison of Euclidean Gaussian capacities \cite[(1.13)]{LTU} and the identity $\capacity(E)=(4\pi)^{d/2}\mathcal C_{1/4}(E)$.
Let $E^\ast=\{(-t,x):(t,x)\in E\}$.
By a change of variables and the time-reversal invariance of thermal capacity \cite[Lemma~2.1]{EG1982},
\[
\mathcal C_a^\dagger(E)=\mathcal C_a(E^\ast)\asymp\capacity(E^\ast)=\capacity(E).
\]

Monotonicity follows directly from the definition.
For countable subadditivity, suppose $E\subset\bigcup_j E_j$, where $E$ and $E_j$ are compact.
Let $\mu$ be a nonnegative measure supported on $E$ with $G_{1/4}\mu \le 1$.
Partition $E$ into disjoint Borel sets $A_j\subset E_j$ and set $\mu_j=\mu\vert_{A_j}$.
Since $\supp\mu_j\subset E_j$ and $G_{1/4}\mu_j \le 1$,
\[
\mu(E)=\sum_j\mu_j(E_j)\le\sum_j\mathcal C_{1/4}(E_j).
\]
Taking the supremum over $\mu$ proves subadditivity.

Let $E_m\downarrow E$ be decreasing compact sets and choose nonnegative measures $\mu_m$ supported on $E_m$ such that
\[
G_{1/4}\mu_m \le 1,\qquad \mu_m(E_m)\ge\mathcal C_{1/4}(E_m)-1/m.
\]
Passing to a subsequence, $\mu_m$ converges weakly to a nonnegative measure $\mu$ supported on $E$.
By lower semicontinuity of $G_{1/4}$, we have $G_{1/4}\mu \le 1$, and hence
\[
\lim_m\mathcal C_{1/4}(E_m)=\mu(E)\le\mathcal C_{1/4}(E).
\]
The reverse inequality follows from monotonicity.

The scaling identity \eqref{eq1717mon} follows from $G_{1/4}(r^2t,rx)=r^{-d}G_{1/4}(t,x)$.
Finally, no positive multiple of a Dirac measure satisfies $G_{1/4}\mu\le 1$, so every singleton has zero thermal capacity.
\end{proof}

\begin{proof}[Proof of Lemma~\ref{lem02}]
Fix $z\in\mathbb R^{d+1}$ and $\epsilon>0$.
Choose $P$-superparabolic functions $v_j$ with $0\le v_j \le 1$ such that $v_j=1$ on $E_j$ and $v_j(z) \le \red_{E_j}(z)+\epsilon2^{-j}$.
By Lemma~\ref{lem03},
\[
w=\lim_{n\to\infty}\min\left\{1,\sum_{j=1}^n v_j\right\}
\]
is $P$-superparabolic.
Since $w=1$ on $E$, we have
\[
\red_E(z)\le w(z)\le\sum_{j=1}^\infty\red_{E_j}(z)+\epsilon.
\]
Letting $\epsilon\downarrow0$ proves the first assertion.

For the second assertion, let $v$ be nonnegative and $P$-superparabolic in $\mathbb R^{d+1}$ with $v \ge 1$ on $F$.
For $\epsilon\in(0,1)$, the set $\{v>1-\epsilon\}$ is an open neighborhood of $F$ and hence contains $K_m$ for all sufficiently large $m$.
Consequently,
\[
\red_F\le\lim_m\red_{K_m}\le\frac{v}{1-\epsilon}.
\]
Taking the infimum over $v$ and letting $\epsilon\downarrow 0$ proves the assertion.
\end{proof}

\begin{proof}[Proof of Lemma~\ref{lem05}]
Let $\mathcal F$ be the family of $P$-superparabolic functions $v$ in $\mathbb R^{d+1}$ with $0\le v\le 1$ and $v=1$ on $K$.
By truncation at $1$, we have $\red_K=\inf_{v\in\mathcal F}v$.
The first three assertions follow immediately.

To show that $\bal_K$ is $P$-superparabolic, let $h\in C(\overline Q)$ be $P$-caloric in a cylinder $Q$ with $h \le \bal_K$ on $\partial_p Q$.
Since $\bal_K\le v$ for every $v\in\mathcal F$, comparison gives $h\le\red_K$ in $Q$.
By continuity of $h$ and the definition of $\bal_K$, we obtain $h\le\bal_K$ in $Q$.

Finally, by Lemma~\ref{lem03}, $\mathcal F$ is closed under finite minima and Poisson modifications in cylinders
$Q\Subset K^\complement$.
Lemma~\ref{lem04} implies that $\red_K$ is $P$-caloric in $K^\complement$.
In particular, $\red_K$ is continuous there, and hence $\red_K=\bal_K$ in $K^\complement$.
\end{proof}


\end{document}